\documentclass[12pt]{amsart}
\usepackage{amsmath,amssymb,times,color}
\usepackage{graphicx}
\usepackage{dsfont}
\definecolor{webred}{rgb}{0.75,0,0}
\definecolor{webgreen}{rgb}{0,0.75,0}
\usepackage[citecolor=webgreen,colorlinks=true,linkcolor=webred]{hyperref}

\oddsidemargin=\evensidemargin

\newtheorem{thm}{Theorem}[section]
\newtheorem{cor}[thm]{Corollary}
\newtheorem{lem}[thm]{Lemma}
\newtheorem{prop}[thm]{Proposition}

\newtheorem{hyp}[thm]{Hypothesis}

\theoremstyle{definition}

\theoremstyle{remark}
\newtheorem{rem}{Remark}[section]

\numberwithin{equation}{section}

\newcommand{\grad}{\operatorname{\mathbf{grad}}}

\newcommand{\Div}{\operatorname{\mathrm{div}}}

\newcommand{\rot}{\operatorname{\mathrm{curl}}}

\newcommand{\db}[1]{_{\raise-0.3ex\hbox{$\scriptstyle #1$}}}
\newcommand{\dd}[1]{_{\raise-1.5pt\hbox{$\scriptstyle #1$}}}

\newcommand{\di}{\displaystyle}

\newcommand{\dr}{{\rm d}}

\newcommand{\cd}{^{\sf -}}
\newcommand{\iso}{_{\sf +}}
\newcommand{\con}{_{\sf -}}
\newcommand{\is}{^{\sf +}}

\newcommand  {\C}{{\mathbb C}}
\newcommand  {\N}{{\mathbb N}}

\newcommand  {\R}{{\mathbb R}}

\newcommand  {\CC}{\boldsymbol{\mathsf C}}
\newcommand   {\D}{\boldsymbol{\mathsf D}}
\newcommand  {\DD}{{\mathsf D}}
\newcommand  {\EE}{\boldsymbol{\mathsf E}}

\renewcommand  {\H}{{\mathrm H}}
\newcommand  {\HH}{\boldsymbol{\mathsf H}}

\renewcommand  {\L}{{\mathrm L}}

\newcommand  {\NN}{{\mathsf N}}

\renewcommand{\SS}{\boldsymbol{\mathsf S}}
\newcommand  {\RR}{\boldsymbol{\mathsf R}}
\newcommand  {\TT}{{\mathsf T}}

\newcommand  {\WW}{\boldsymbol{\mathsf W}}
\newcommand  {\XX}{\boldsymbol{\mathsf X}}

\newcommand  {\g}{{\boldsymbol{\mathsf g}}}

\newcommand  {\jj}{{\boldsymbol{\mathsf j}}}

\newcommand  {\nn}{{\boldsymbol{\mathsf n}}}

\newcommand  {\uu}{\boldsymbol{\mathsf u}}

\newcommand  {\ww}{\boldsymbol{\mathsf w}}
\newcommand  {\xx}{\boldsymbol{\mathsf x}}
\newcommand  {\yy}{\boldsymbol{\mathsf y}}

\newcommand  {\jjw}{\widetilde{\boldsymbol{\mathsf\jmath}}}

\newcommand  {\cU}{{\mathcal U}}

\newcommand  {\ba}{\underline a}
\newcommand  {\bal}{\underline\alpha}
\newcommand  {\bs}{\underline \sigma}
\newcommand  {\bH}{\mathbf{H}}
\newcommand  {\bL}{\mathbf{L}}
\newcommand  {\gR}{\mathfrak{R}}
\newcommand  {\gT}{\mathfrak{T}}
\def\on#1{\!\left.\vphantom{|_|}\right|_{#1}}  
\renewcommand{\Re}{\operatorname{Re}}
\renewcommand{\Im}{\operatorname{Im}}

\begin{document}

\title[Uniform estimates]{Uniform estimates for transmission problems with high contrast in heat conduction and electromagnetism }

\author{Gabriel Caloz, Monique Dauge, Victor Peron}
\address{Labortatoire IRMAR, UMR 6625 du CNRS,
Campus de Beaulieu
35042 Rennes cedex, France}



\begin{abstract}
In this paper we prove uniform a priori estimates for transmission
problems with constant coefficients on two subdomains, with a special emphasis 
for the case when the ratio between these coefficients is large.
In the most part of the work, the interface between the two subdomains is supposed to be Lipschitz.
We first study a scalar transmission problem
which is handled through a converging asymptotic series.
Then we derive uniform a priori estimates for Maxwell transmission
problem set on a domain made up of a dielectric and a 
highly conducting material. The technique is based 
on an appropriate decomposition of the electric field, whose gradient part is
estimated thanks to the first part. As an application, we develop an argument for the convergence
of an asymptotic expansion as the conductivity tends to infinity. 
\end{abstract}

\maketitle

\section{Introduction}
\label{S1}
The goal of our work is to derive uniform a priori estimates for 
transmission problems in media presenting high contrast in their
material properties. We investigate in particular the heat transfer 
equation
\begin{equation}
\label{0E1}
   \Div \ba \grad \varphi = f
\end{equation}
and the Maxwell equations given by Faraday's and Amp\`ere's laws
\begin{equation}
\label{0E2}
    \rot  \EE - i \omega\mu_0 \HH = 0 \quad \mbox{and}\quad
    \rot \HH + (i\omega\varepsilon_0 - \bs) \EE =   \jj\,.
\end{equation}
Here, $\ba$ represents the heat conductivity and $\bs$ the electrical conductivity.
We assume that these equations are set in a domain $\Omega$ made up of two subdomains
$\Omega_+$ and $\Omega_-$ in which the coefficients $\ba$ and $\bs$ take two different values $(a_+,a_-)$ and $(\sigma_+,\sigma_-)$, respectively. These equations are complemented by suitable boundary conditions. 
Our interest is their solvability together with uniform energy or regularity estimates, namely
\begin{itemize}
\item when the ratio $|a_-|/|a_+|$ tends to infinity in the case of eq.\ \eqref{0E1}
\item when $\sigma_+=0$ (insulating or dielectric material) and $\sigma_-\equiv\sigma$ tends to infinity (highly conducting material) in the case of eq.\ \eqref{0E2}.
\end{itemize}

We address different, though connected, issues for these two problems, namely the issue of {\em uniform piecewise regularity} in Sobolev norms for solutions $\varphi$ of equation \eqref{0E1}, and the issue of {\em uniform $\L^2$ estimates} for the electromagnetic field $(\EE,\HH)$ solution of system \eqref{0E2}. None of these questions have obvious answers, all the more since we do not assume that the interface $\Sigma$ between $\Omega_+$ and $\Omega_-$ is smooth.

Our whole analysis is valid under the only following assumption on the interface
\begin{equation}
\label{0E3}
   \Sigma\quad \mbox{is a bounded Lipschitz surface in $\R^3$.}
\end{equation}
In the Maxwell case, similar estimates as ours are obtained in \cite{H-J-N07}, but under a stronger regularity assumption on $\Sigma$. Our approach differs, being based on a decomposition of
the electric field given in \cite{ABDG98}. The gradient part of the decomposition
is handled through the uniform regularity estimates proved for equation \eqref{0E1}. 

The paper is organized as follows. In section \ref{S12} we introduce the notations and 
give the main results. In section \ref{S3} we prove uniform piecewise $\H^{\frac32}$ estimates for solutions of the scalar interface problem \eqref{0E1} with exterior Dirichlet or Neumann boundary 
conditions. In section \ref{S4} we prove uniform estimates for the electromagnetic field $(\EE,\HH)$ solution of the Maxwell system \eqref{0E2} when the conductivity $\sigma$ of the conducting part is high. We conclude our paper in section \ref{S5} by an application of the previous uniform estimates to the convergence study of an asymptotic expansion as the conductivity tends to infinity.

\section{Notations and main results}
\label{S12}

\begin{figure}[h]
\includegraphics[scale=1.0]{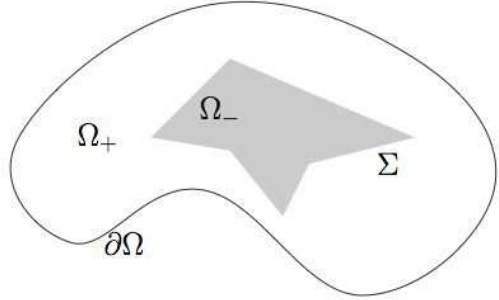}
\vskip -3mm
\caption{The domain $\Omega$ and its subdomains $\Omega_+$ and $\Omega_-$}
\label{F1}
\end{figure}

Let $\Omega$ be a smooth bounded simply connected domain in ${\R}^3$ with boundary $\partial\Omega$, and $\Omega_{-}\subset\subset\Omega$ be a Lipschitz connected subdomain of $\Omega$, with boundary $\Sigma$. We denote by $\Omega_{+}$ the complementary of $\overline{\Omega}_{-}$ in $\Omega$, see Figure \ref{F1}.

We denote by $\varphi^+$ ($\varphi^-$) the restriction of any function $\varphi$ to $\Omega_{+}$ ($\Omega_{-}$).

\subsection{Scalar problem}
We consider both Dirichlet and Neumann external boundary conditions associated with equation \eqref{0E1} and introduce the functional spaces suitable for their variational formulation: $V_\DD=\H^1_0(\Omega)$ for Dirichlet and 
$V_\NN=\{ \varphi \in  \H^1(\Omega)\, | \, \int_\Omega  \varphi \, \dr\xx= 0\}$ for Neumann. For any given function $\ba=(a_+,a_-)$ determined by the two constants $a_\pm$ on $\Omega_\pm$ and in either case ($V=V_{\NN}$ or $V=V_{\DD}$) the variational problem is: Find $\varphi \in V$ such that
\begin{multline}
   \forall \psi \in  V, \quad 
   \int_{\Omega_{+}} a_+\nabla\varphi^+ \cdot \nabla\overline{\psi}{}^+ \,\dr\xx
   + \int_{\Omega_{-}} a_-\nabla\varphi^- \cdot \nabla\overline{\psi}{}^-\, \dr\xx= \\[-1ex]
   -\int_{\Omega}  f \; \overline{\psi}\,\dr\xx 
   + (a_+-a_-) \!\int_{\Sigma} g\; \overline{\psi} \, \dr s,
\label{1E1}
\end{multline}
where the right-hand side $(f,g)$ satisfies the regularity assumption
\begin{equation}
\label{1E7}
f\in\L^2(\Omega)\quad \mbox{and} \quad g\in \L^2(\Sigma)
\end{equation}
and the extra compatibility conditions
\begin{equation}
 \int_{\Omega}  f\, \dr\xx =0 \ \ \ \mbox{and} \ 
 \int_{\Sigma} g \, \dr s = 0\ \mbox{ if } \ V=V_\NN \quad \mbox{and} \quad 
 \int_{\Sigma} g \, \dr s =0 \ \mbox{ if } \ V= V_\DD\,.
\label{1E2}
\end{equation}
Our main result in the scalar case is the following piecewise $\H^{3/2}$ a priori estimate, uniform with respect to the ratio $\rho:=a_-(a_+)^{-1}$.
It applies both to Dirichlet and Neumann boundary conditions.

\begin{thm}
\label{2T9}
Let us assume that $a_+\neq0$.
There exist a constant $\rho_{0}>0$ independent of $a_+$ such that for all $a_-\in\{z\in\C \, | \,|z|\geqslant \rho_{0}|a_+| \}$, the problem \eqref{1E1} with data $(f,g)$ satisfying \eqref{1E7}-\eqref{1E2} has a unique solution $\varphi\in V$, which moreover is piecewise $\H^{3/2}$ and satisfies the uniform estimate
\begin{equation}
\label{ue}
   \|\varphi^+\|_{\frac32,\Omega_+} + \|\varphi^-\|_{\frac32,\Omega_-}
   \leqslant C_{\rho_{0}} \big( |a_+|^{-1} \|f \|_{0,\Omega} + \|g\|_{0,\Sigma}\big)
   \end{equation}
with a constant $C_{\rho_{0}}>0$, independent of $a_+$, $a_-$, $f$, and $g$.
\end{thm}

This statement is proved in the next section using an asymptotic expansion for $\varphi$ with respect to the powers of $\rho^{-1}=a_+(a_-)^{-1}$. The estimate \eqref{ue} will be a consequence of the {\em convergence} of this series in the piecewise $\H^{3/2}$-norm. The dependence of $\rho_0$ and $C_{\rho_0}$ on the overall configuration is discussed in Remark \ref{3R1} after the proof.

\begin{rem}
\label{2R1}
The estimate \eqref{ue} is uniform for fixed $a_+$ when $|a_-|$ tends to infinity. The roles of $a_+$ and $a_-$ can be exchanged and an estimate similar to \eqref{ue} proved. In fact, there holds a more precise estimate where $a_+$ and $a_-$ play symmetric roles, see Proposition \ref{3P1}.
\end{rem}

\begin{rem}
\label{2R1b}
In the Neumann case, the compatibility conditions \eqref{1E2} are necessary for the right hand side of problem \eqref{1E1} to be compatible for all values of $(a_-,a_+)$, because of the factor $(a_+-a_-)$ in front of the integral on $\Sigma$. If this factor is replaced by $1$, then, under the weaker conditions
\begin{equation}
- \int_{\Omega}  f\, \dr\xx +
 \int_{\Sigma} g \, \dr s = 0\ \mbox{ if } \ V=V_\NN \quad \mbox{and} \quad 
 \mbox{nothing} \ \mbox{ if } \ V= V_\DD\,,
\label{1E2b}
\end{equation}
the problem \eqref{1E1} is still solvable for $\rho$ large enough, see Proposition \ref{3P2}.
\end{rem}

\begin{figure}[h]
\includegraphics[scale=1.0]{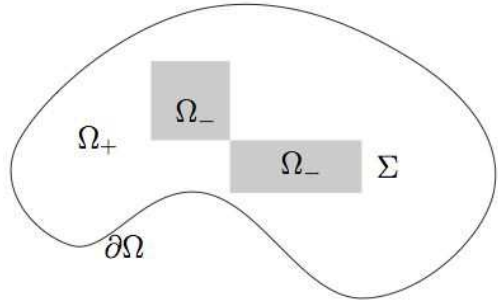}
\vskip -3mm
\caption{A non-Lipschitz interface $\Sigma$}
\label{F2}
\end{figure}

\begin{rem}
\label{2R2}
The assumption that $\Sigma$ is Lipschitz is necessary: There exist non-Lipschitz interfaces such that estimate \eqref{ue} does not hold. In two dimensions of space, such an example is provided by the checkerboard configuration (Figure \ref{F2}), cf.\ \cite[Theorem 8.1]{C-D-N99}.
\end{rem}

\begin{rem}
\label{2R3}
If the Lipschitz interface $\Sigma$ is {\em polyhedral}, there hold uniform piecewise $\H^s$ estimates for any exponent $s$, $\frac 32<s\leqslant s_\Sigma$, with some $\frac32<s_\Sigma\leqslant2$, cf.\ \cite[Ch.\,1, \S 1.5]{Peron09}:
\begin{equation}
\label{ues}
   \|\varphi^+\|_{s,\Omega_+} + \|\varphi^-\|_{s,\Omega_-}
   \leqslant C_{\rho_{0},s} 
   \big(|a_+|^{-1} \|f \|_{s-2,\Omega} + \|g\|_{s-\frac32,\Sigma}\big).
\end{equation}
This is a consequence of our proof (see Remark \ref{3R1}) in combination with elliptic estimates in polyhedral domain, cf.\ \cite{Dauge88}. In particular, if $\Sigma$ contains some edge, then $s_\Sigma<2$.
\end{rem}

\begin{rem}
\label{2R4}
If the interface $\Sigma$ is {\em smooth}, there hold uniform piecewise $\H^s$ estimates  for any $s\geqslant2$:
\begin{equation}
\label{uesb}
   \|\varphi^+\|_{s,\Omega_+} + \|\varphi^-\|_{s,\Omega_-}
   \leqslant C_{\rho_{0},s} 
   \big(|a_+|^{-1} \|f \|_{s-2,\Omega} + \|g\|_{s-\frac32,\Sigma}\big).
\end{equation}
This is again a consequence of our proof in combination with standard elliptic estimates, cf.\ \cite{Agmon1965} for instance.
\end{rem}

\subsection{Maxwell problem}
We consider two types of boundary conditions to complement the Maxwell harmonic equations \eqref{0E2} on $\partial\Omega$: Either the perfectly insulating conditions
\begin{subequations}
\begin{equation}
\label{ME2}
   \EE\cdot\nn=0 \quad\mbox{and}\quad \HH\times\nn=0\quad\mbox{on}\quad\partial\Omega\,,
\end{equation}
where $\nn$ denotes the outer normal vector, or the perfectly conducting conditions
\begin{equation}
\label{ME1}
   \EE\times\nn=0 \quad\mbox{and}\quad \HH\cdot\nn=0\quad\mbox{on}\quad\partial\Omega\,.
\end{equation}
\end{subequations}
In both cases, for the conductivity $\bs=(0,\sigma)$, we can prove uniform a priori estimates for the electromagnetic field as $\sigma\to\infty$ provided the following condition on limit problems in the dielectric part $\Omega_+$ is valid: 

\begin{hyp}
\label{H1}
The angular frequency $\omega$ is not an eigenfrequency of the problem
\begin{equation}
\label{PH1}
 \left\{
   \begin{array}{lll}
    \rot  \EE - i \omega\mu_0 \HH = 0 \quad \mbox{and}\quad
    \rot \HH + i\omega\varepsilon_0 \EE = 0 \quad&\mbox{in}\quad \Omega_{+}
\\[5pt]
\EE\times \nn  = 0 \quad \mbox{and} \quad
\HH\cdot\nn=0  \quad &\mbox{on}\quad \Sigma
\\[5pt]
\eqref{ME2}\mbox{ or }\eqref{ME1}  \quad &\mbox{on}\quad \partial\Omega.
   \end{array}
    \right.
\end{equation}
\end{hyp}

Our main result for Maxwell equations is the following 	a priori estimate, uniform as $\sigma\to\infty$. The right hand side $\jj$ is chosen to belong in $\bH(\Div,\Omega)$ or $\bH_0(\Div,\Omega)$ where:
\begin{subequations}
\begin{gather}
   \bH(\Div,\Omega) = \{\uu\in\bL^2(\Omega)\; |\quad \Div\uu\in\L^2(\Omega) \} 
   \quad\mbox{with}\quad \bL^2(\Omega) = \L^2(\Omega)^3\\
   \bH_0(\Div,\Omega) = \{\uu\in\bH(\Div,\Omega)\; |\quad 
   \uu\cdot\nn=0 \ \ \mbox{ on } \ \partial\Omega \}.
\end{gather}
\end{subequations}

\begin{thm}
\label{4T0}
Under Hypothesis {\em\ref{H1}}, there are constants $\sigma_{0}$ and $C>0$, such that for all $\sigma\geqslant\sigma_0$, the Maxwell problem \eqref{0E2} with boundary condition \eqref{ME2} and data $\jj\in\bH_0(\Div,\Omega)$ has a unique solution $(\EE,\HH)$ in $\bL^2(\Omega)^2$, which satisfies:
\begin{equation}
   \|\EE\|_{0,\Omega} + \|\HH\|_{0,\Omega} 
   + \sqrt\sigma\, \|\EE\|_{0,\Omega_-} \leqslant C \|\jj\|_{\bH(\Div,\Omega)}.
\label{3E1}
\end{equation}
A similar result holds for boundary conditions \eqref{ME1} and data $\jj\in\bH(\Div,\Omega)$.
\end{thm}

This theorem is proved in section \ref{S4}. It is based in particular on a decomposition of the electric field into a regular field $\ww$ in $\bH^1(\Omega)$ and a gradient field $\nabla\varphi$ for which Theorem \ref{2T9} will be used.

\section{Proof of uniform scalar regularity estimates}
\label{S3}

\subsection{The problem}
\label{S3.0}
Under the assumptions of Theorem \ref{2T9}, we normalize the equations by dividing by $a_+$. Denoting the quotient $a_-(a_+)^{-1}$ by $\rho$, and still denoting by $f$  the new right hand side (i.e., the old one divided by $a_+$), we can write problem \eqref{1E1} in the form of the following transmission problem
\begin{equation}
\label{PFort}
\left\{
   \begin{array}{lll}
 \rho\Delta \varphi^{-}    &= \quad f^{-} \quad
&\mbox{in}\quad \Omega_{-}
\\[0.5ex]
 \Delta \varphi^{+}   &=\quad f^{+} \quad
&\mbox{in}\quad \Omega_{+}
\\[0.5ex]
\varphi^{-}    & =\quad \varphi^{+}  \quad &\mbox{on}\quad \Sigma
\\[0.5ex]
\partial_{\nn}\varphi^+ - \rho\partial_{\nn}\varphi^-  & =\quad (1-\rho) g \quad
&\mbox{on}\quad \Sigma
\\[0.5ex]
 \mbox{(b.c.)} &\quad
&\mbox{on}\quad  \partial\Omega
   \end{array}
    \right.
\end{equation}
where $\partial_{\nn}$ denotes the normal derivative (inner for $ \Omega_{-}$, outer for $ \Omega_{+}$).
The external boundary conditions (b.c.) are either Neumann or Dirichlet conditions.

Our method of proof for Theorem \ref{2T9} consists in the determination of a series expansion 
in powers of $\rho^{-1}$ for $\varphi$ solution of \eqref{PFort}: We are looking for solutions in the form of power series
\begin{equation}
\varphi=
\left\{
   \begin{array}{lll}
   \di\sum_{n=0}^\infty\, \varphi_{n}^{+} \,\rho^{-n}\quad
   &\mbox{in}\quad \Omega_{+}
\\[0.5ex]
   \di\sum_{n=0}^\infty\,\varphi_{n}^{-}\, \rho^{-n}\quad
   &\mbox{in}\quad \Omega_{-}\,.
  \end{array}
    \right.
\label{1E6}
\end{equation}
Since the expansions are different according to the boundary conditions, we treat first the Neumann case in subsection \ref{S3.1} and the Dirichlet case in subsection \ref{S3.2}. We prove complementary results in subsection \ref{S3.3}.

\subsection{Neumann external b.c.}
\label{S3.1}
Inserting the ansatz \eqref{1E6}
in the system \eqref{PFort}, we get the following families of problems, coupled by their conditions on $\Sigma$:
\begin{equation}
\left\{
   \begin{array}{lll}
\Delta \varphi_{0}^-& =\quad0 \quad
&\mbox{in}\quad \Omega_{-}
\\[0.5ex]
\partial_{\nn}\varphi_{0}^-& =\quad g \quad
&\mbox{on}\quad \Sigma
   \end{array}
    \right.
 \label{1E40}
\end{equation}
and
\begin{equation}
\left\{
   \begin{array}{lll}
\Delta \varphi_{0}^+ & = \quad f^{+} \quad
&\mbox{in}\quad \Omega_{+}
\\[0.5ex]
\varphi_{0}^+ & = \quad\varphi_{0}^-  \quad &\mbox{on}\quad \Sigma
\\[0.5ex]
\partial_{\nn}\varphi_{0}^+&=\quad 0 \quad
&\mbox{on}\quad \partial\Omega
   \end{array}
    \right.
 \label{1E41}
\end{equation}
and for $k\in\N^*$ (here $\delta_{k}^l$ is the Kronecker symbol)
\begin{equation}
\left\{
   \begin{array}{lll}
   \Delta \varphi_{k}^-&=\quad\delta_{k}^1\,f^{-} \quad
   &\mbox{in}\quad \Omega_{-}
\\[0.5ex]
   \partial_{\nn}\varphi_{k}^-&=\quad - \delta_{k}^1\,g + \partial_{\nn}\varphi_{k-1}^+
   &\mbox{on}\quad \Sigma
   \end{array}
    \right.
 \label{1E42}
\end{equation}
and
\begin{equation}
\left\{
   \begin{array}{lll}
\Delta \varphi_{k}^+ & =\quad 0 \quad
&\mbox{in}\quad \Omega_{+}
\\[0.5ex]
\varphi_{k}^+  &= \quad \varphi_{k}^-  \quad &\mbox{on}\quad \Sigma
\\[0.5ex]
\partial_{\nn}\varphi_{k}^+&=\quad 0 \quad
&\mbox{on}\quad \partial\Omega\,.
   \end{array}
    \right.
 \label{1E43}
\end{equation}
Thus we alternate the solution of a Neumann problem in $\Omega_-$ and a mixed Dirichlet-Neumann problem in $\Omega_+$. Since we have assumed that $\Sigma$ is a Lipschitz surface, we have a precise and optimal functional framework to describe these operators and their inverse.

We need the notation ($s<2$)
\begin{equation}
\label{1E5}
   \H^s(U,\Delta) := \{\varphi\in\H^s(U)\;|\;\; \Delta \varphi\in\L^2(U)\}.
\end{equation}
Using \cite[cor.\ 5.7]{JerisonKenig82} (see also \cite{CoDa98} for a similar context) with the fact that $\Sigma$ is Lipschitz, we obtain the following equalities between spaces on $\Omega_-$:
\begin{subequations}
\label{3E8}
\begin{eqnarray}
   \big\{ \varphi\in \H^1(\Omega_-,\Delta),\;
   \partial_\nn\varphi\on{\Sigma}\in \L^2(\Sigma) \big\} & = &
   \H^{\frac32}(\Omega_-,\Delta)\\
   \big\{ \varphi\in \H^1(\Omega_-,\Delta),\;
   \varphi\on{\Sigma}\in \H^1(\Sigma) \big\} & = &
   \H^{\frac32}(\Omega_-,\Delta)
\end{eqnarray}
\end{subequations}
and on $\Omega_+$
\begin{subequations}
\label{3E9}
\begin{eqnarray}
   \qquad\big\{ \varphi\in \H^1(\Omega_+,\Delta),\;
   \varphi\on{\Sigma}\in \H^1(\Sigma), \; 
   \partial_\nn\varphi\on{\partial\Omega}\in \L^2(\partial\Omega) \big\} & = &
   \H^{\frac32}(\Omega_+,\Delta)\\
   \qquad\big\{ \varphi\in \H^1(\Omega_+,\Delta),\;
   \partial_\nn\varphi\on{\Sigma}\in \L^2(\Sigma), \; 
   \varphi\on{\partial\Omega}\in \H^1(\partial\Omega) \big\} & = &
   \H^{\frac32}(\Omega_+,\Delta).   
\end{eqnarray}
\end{subequations}

As a consequence of the previous equalities, the following definitions for the resolvent operators $\gR_-$ and $\gR_+$ make sense and define bounded operators: \\
$\bullet$ \ $\gR_-$ is the resolvent of the Neumann problem on $\Omega_-$
\begin{multline}
\label{Rm}
   \gR_-: \Big\{(F,G)\in\L^2(\Omega_-)\times\L^2(\Sigma)\;|\;\;
   \int_{\Omega_-}\!F\,\dr\xx + \int_\Sigma G\,\dr s=0\Big\} \\
   \longrightarrow\quad
   \Big\{\Phi\in\H^{\frac32}(\Omega_-,\Delta)\ |\;\;\int_{\Omega_-}\!\Phi\,\dr\xx=0\Big\}
\end{multline}
where $\Phi=\gR_-(F,G)$ satisfies $\Delta\Phi=F$ in $\Omega_-$ and $\partial_\nn\Phi=G$ on $\Sigma$,
and\\
$\bullet$ \ $\gR_+$ is the resolvent of the Dirichlet-Neumann problem on $\Omega_+$
\begin{equation}
\label{Rp}\textstyle
   \gR_+: \Big\{(F,H)\in\L^2(\Omega_+)\times\H^1(\Sigma)\Big\} \quad\longrightarrow\quad
   \H^{\frac32}(\Omega_+,\Delta)
\end{equation}
where $\Phi=\gR_+(F,H)$ satisfies $\Delta\Phi=F$ in $\Omega_+$, $\Phi=H$ on $\Sigma$ and $\partial_\nn\Phi=0$ on $\partial\Omega$.

A further consequence of equalities \eqref{3E8}-\eqref{3E9} is that the following trace operators make sense and define bounded operators: \\ 
$\bullet$ \ $\gT^0_-$ is the Dirichlet trace on $\Sigma$ from inside $\Omega_-$: $\varphi^-\mapsto\varphi^-|_\Sigma$ and is bounded
\begin{equation}
\label{Tm}
   \gT^0_-:\H^{\frac32}(\Omega_-,\Delta) \quad \longrightarrow\quad \H^1(\Sigma)
\end{equation}
$\bullet$ \ $\gT^1_+$ is the Neumann trace on $\Sigma$ from inside $\Omega_+$: $\varphi^+\mapsto\partial_\nn\varphi^+|_\Sigma$ and is bounded
\begin{equation}
\label{Tp}
   \gT^1_+:\H^{\frac32}(\Omega_+,\Delta) \quad \longrightarrow\quad \L^2(\Sigma)
\end{equation}
Note that none of the two operators $\gT^0_-$ or $\gT^1_+$ is bounded if $\H^{\frac32}(\Omega_\pm,\Delta)$ is replaced by the larger space $\H^{\frac32}(\Omega_\pm)$.

\subsubsection{Discussion of elementary problems.}
\label{S321}
The Neumann problem \eqref{1E40} admits the solution $\varphi_{0}^-=\gR_-(0,g)$, since the compatibility condition $\int_{\Sigma} g \; \dr s = 0$ holds by assumption. We have the estimate
\begin{equation}
\label{1E57a}
   \|\varphi_{0}^-\|_{\H^{\frac32}(\Omega_-,\Delta)} \leqslant  
   C_{-} \|g\|_{0,\Sigma} \,,
\end{equation}
where $C_{-}$ is the operator norm of the operator $\gR_-$.

Thus $\varphi_{0}^-$ belongs to $\H^{\frac32}(\Omega_-,\Delta)$. Hence its trace $\gT^0_-\varphi^-_0=\varphi^-_0|_\Sigma$ belongs to $\H^1(\Sigma)$.
Next, problem \eqref{1E41} admits the solution $\varphi_{0}^+ = \gR_+(f^+,\varphi^-_0|_\Sigma)$, and $\varphi^+_0$ belongs to $\H^{\frac32}(\Omega_+,\Delta)$ with the estimates
\begin{equation}
\label{1E59a}
   \|{\varphi}_{0}^+\|_{\H^{\frac32}(\Omega_+,\Delta)} \leqslant  
   C_{+} \big( \|f^{+}\|_{0,\Omega_{+}} + \|{\varphi}_{0}^-\|_{1,\Sigma} \big).
\end{equation}
Here $C_{+}$ is the operator norm of the operator $\gR_+$.

Then we continue with problems \eqref{1E42} and \eqref{1E43} in a similar way, the only point to discuss being the compatibility condition in the Neumann problem \eqref{1E42}.

\begin{lem}
The Neumann problem \eqref{1E42} is compatible.
\end{lem}

\begin{proof}
For $k=1$, we must show that
\begin{equation}
\label{1e54}
   \int_{\Omega_{-}}f^{-}\dr\xx+\int_{\Sigma} (-g+\partial_{\nn}\varphi_{0}^+)\, \dr s=0\,.
\end{equation}
According to \eqref{1E41}, $\Delta\varphi_{0}^+=f_{+}$ in $\Omega_+$ and $\partial_{\nn}\varphi_{0}^+=0$ on $\partial\Omega$.
Integrating by parts, we get
\begin{equation*}
   \int_{\Sigma} \partial_{\nn}\varphi_{0}^+\, \dr s=\int_{\Omega_{+}}f_{+}\dr\xx\,.
\end{equation*}
Thus, we deduce from hypothesis \eqref{1E2} the compatibility condition \eqref{1e54}.

For $k\geqslant2$, let us assume that the term $\varphi_{k-1}^+$ was built. We must show that
\begin{equation}
\label{CCompa}
   \int_{\Sigma} \partial_{\nn}\varphi_{k-1}^+\, \dr s=0\,.
\end{equation}
According to \eqref{1E43}, $\Delta\varphi_{k-1}^+=0$ in $\Omega_+$ and $\partial_{\nn}\varphi_{k-1}^+=0$ on $\partial\Omega$.
Integrating by parts in $\Omega_{+}$, we get \eqref{CCompa}.
\end{proof}

Consequently, for all $k\geqslant 1$, the Neumann problem \eqref{1E42} admits the solution $\varphi_{k}^-:=\gR_-(\delta_{k}^1 f^{-}, - \delta_{k}^1 g+\gT^1_+\varphi_{k-1}^+)$. Then $\varphi_{k}^-$ belongs to $\H^{\frac32}(\Omega_-,\Delta)$ with the estimates
\begin{equation}
\label{1E57}
   \|\varphi_{k}^-\|_{\H^{\frac32}(\Omega_-,\Delta)} \leqslant  
   C_{-} \Big(\delta_{k}^1 \big(\|f^{-}\|_{0,\Omega_{-}}+\|g\|_{0,\Sigma} \big)
   +\|\partial_{\nn}\varphi_{k-1}^+\|_{0,\Sigma}\Big) .
\end{equation}
Finally problem \eqref{1E43} admits a unique solution $\varphi_{k}^+\in\H^{\frac32}(\Omega_+,\Delta)$ with the estimate
\begin{equation}
\label{1E59}
   \|{\varphi}_{k}^+\|_{\H^{\frac32}(\Omega_+,\Delta)} \leqslant  
   C_{+} \|{\varphi}_{k}^-\|_{1,\Sigma} \,.
\end{equation}
In \eqref{1E57} and \eqref{1E59}, the constants $C_-$ and $C_+$ are the same as in \eqref{1E57a} and \eqref{1E59a}.

\subsubsection{Uniform estimates}
Let $C^0$ and $C^1>0$ be the operator norms of the Dirichlet and Neumann traces $\gT^0_-$ and $\gT^1_+$, respectively cf.\ \eqref{Tm} and \eqref{Tp}. 
We set $\alpha= C_+C^0C_{-}C^{1}$, with the constants $C_{+}$ in \eqref{1E59a} and $C_{-}$ in \eqref{1E57a}. According to \eqref{1E57} and \eqref{1E59}, we see by an induction on $n\in\N^*$ that
\begin{equation}
\label{eq:3.20}
\left\{
   \begin{array}{lll}
   \|\varphi_{n}^- \|_{\H^{\frac32}(\Omega_-,\Delta)} \leqslant 
   \alpha^{n-1}\|\varphi_{1}^- \|_{\H^{\frac32}(\Omega_-,\Delta)}
\\[1.5ex]
   \|\varphi_{n}^+ \|_{\H^{\frac32}(\Omega_+,\Delta)}\leqslant 
   C_{+}C^0 \alpha^{n-1}\|\varphi_{1}^- \|_{\H^{\frac32}(\Omega_-,\Delta)}.
   \end{array}
    \right.
\end{equation}
Let $\rho_{0}>0$ such that $\rho_{0}^{-1}\alpha<1$. Then, for all $\rho\in\C$ such that $|\rho|\geqslant\rho_{0}$, the series of general terms $\rho^{-n}\varphi_{n}^{-}$ et $\rho^{-n}\varphi_{n}^{+}$ converge respectively in $\H^{\frac32}(\Omega_-,\Delta)$ and $\H^{\frac32}(\Omega_+,\Delta)$. We denote by $\varphi_{(\rho)}$ the sum of these series. Moreover, normal convergence is geometric with common ratio $|\rho^{-1}|\alpha$, bounded by $\rho_{0}^{-1}\alpha$. Hence
\begin{equation}
\left\{
   \begin{array}{lll}
   \|\varphi^-_{(\rho)} \|_{\H^{\frac32}(\Omega_-,\Delta)}\leqslant 
   \frac{\rho_0}{\rho_0 - \alpha}\;|\rho^{-1}|\,
   \|\varphi_{1}^- \|_{\H^{\frac32}(\Omega_-,\Delta)} + 
   \|\varphi_{0}^- \|_{\H^{\frac32}(\Omega_-,\Delta)}\,,
\\[1.5ex]
   \|\varphi^+_{(\rho)} \|_{\H^{\frac32}(\Omega_+,\Delta)} \leqslant 
   C_{+}C^0 \frac{\rho_0}{\rho_0 - \alpha}\;|\rho^{-1}|\,
   \|\varphi_{1}^- \|_{\H^{\frac32}(\Omega_-,\Delta)}  +
   \|\varphi_{0}^+ \|_{\H^{\frac32}(\Omega_+,\Delta)}\,.
   \end{array}
    \right.
\label{1E52}
\end{equation}
According to \eqref{1E57} for $k=1$, and \eqref{1E57a}-\eqref{1E59a} for $k=0$,
\begin{equation}
   \|\varphi_{1}^- \|_{\H^{\frac32}(\Omega_-,\Delta)} \leqslant
   C_{-} \Big(\|f^{-}\|_{0,\Omega_{-}}+\|g\|_{0,\Sigma}
   + C^1 C_{+} \big( \|f^{+}\|_{0,\Omega_{+}} + 
   C^0 C_{-} \|g\|_{0,\Sigma} \big) \Big) .
\label{1E55}
\end{equation}
With \eqref{1E57a}, \eqref{1E59a}, \eqref{1E52} and \eqref{1E55}, we deduce the uniform estimate for $|\rho|\geqslant\rho_0$
\begin{equation}
\label{1E56}
   \|\varphi^+_{(\rho)}\|_{\frac32,\Omega_+} + \|\varphi^-_{(\rho)}\|_{\frac32,\Omega_-}
   \leqslant C{(\rho_{0})}\, \big(\|f \|_{0,\Omega} + \|g\|_{0,\Sigma}\big).
\end{equation}

\subsubsection{Proof of Theorem \ref{2T9} in the Neumann case} 
By construction, $\varphi_{(\rho)}$ is solution of the problem \eqref{PFort}. Hence, $\varphi_{(\rho)}\in \H^1(\Omega)$. Setting
\[
   \varphi'_{(\rho)} := \varphi_{(\rho)} - \frac1{|\Omega|} \int_{
   \Omega} \varphi_{(\rho)}(\xx) \, \dr\xx
\]
we obtain a solution $\varphi=\varphi'_{(\rho)}$ of the variational problem \eqref{1E1}, which moreover satisfies estimates \eqref{1E56}, hence estimates \eqref{ue}.

It remains to prove that the solution of the variational problem \eqref{1E1} is {\em unique} when $|\rho|\geqslant\rho_0$. Let $\varphi_*\in V$ be solution of problem \eqref{1E1} for such a $\rho$ and for $f=0$, $g=0$:
\begin{equation}
\label{3E24}
   \forall \psi \in  V, \quad 
   \int_{\Omega_{+}} a_+\nabla\varphi^+_* \cdot \nabla\overline{\psi}{}^+ \,\dr\xx
   + \int_{\Omega_{-}} a_-\nabla\varphi^-_* \cdot \nabla\overline{\psi}{}^-\, \dr\xx
   =   0. 
\end{equation} 
On the other hand, the power series construction yields a solution $\psi_*$ of problem \eqref{1E1} with $\overline a_\pm$ instead of $a_\pm$ and with $f=\varphi_*$, $g=0$ (note that these data satisfy assumption \eqref{1E2}):
\begin{equation}
\label{3E25}
   \forall \varphi \in  V, \quad 
   \int_{\Omega_{+}} \overline a_+\nabla\psi_*^+ \cdot \nabla\overline{\varphi}{}^+ \,\dr\xx
 + \int_{\Omega_{-}} \overline a_-\nabla\psi_*^- \cdot \nabla\overline{\varphi}{}^- \,\dr\xx
 = 
   -\int_{\Omega}  \varphi_* \; \overline{\varphi}\,\dr\xx ,
\end{equation}

Taking the conjugate of \eqref{3E25} for $\varphi=\varphi_*$ and \eqref{3E24} for $\psi=\psi_*$, we find
\[
   \int_{\Omega}  \varphi_* \; \overline{\varphi_*}\,\dr\xx = 0.
\]
Hence the uniqueness, which concludes the proof of Theorem \ref{2T9}.

\begin{rem}
\label{3R1}
From the above proof, we can see that the constants $\rho_0$ and $C_{\rho_0}$ depend only on the four operator norms $C^0$, $C^1$ (trace operators $\gT^0_-$ \eqref{Tm} and $\gT^1_+$  \eqref{Tp}), $C_-$ and $C_+$ (resolvent operators $\gR_-$ \eqref{Rm} and $\gR_+$ \eqref{Rp}). The extension of the estimates \eqref{ue} to different sets of Sobolev indices, cf.\ Remarks \ref{2R3} and \ref{2R4}, depends on the boundedness of the four operators (the orthogonality conditions are understood for the last two ones): 
\[
   \gT^0_-:\H^{s}(\Omega_-) \longrightarrow \H^{s-\frac12}(\Sigma)\, ,
   \quad
   \gT^1_+:\H^{s}(\Omega_+) \longrightarrow \H^{s-\frac32}(\Sigma)\,,
\]
\[
   \gR_-: \H^{s-2}(\Omega_-)\times\H^{s-\frac32}(\Sigma)
   \longrightarrow
   \H^{s}(\Omega_-)\,,\quad
   \gR_+: \H^{s-2}(\Omega_+)\times\H^{s-\frac12}(\Sigma)
   \longrightarrow
   \H^{s}(\Omega_+)\,.
\]
In particular,  none of them is bounded for $s=\frac32$, so we cannot set $s=\frac32$ in estimate \eqref{ues}.
\end{rem}

\subsection{Dirichlet external b.c.}
\label{S3.2}
Here $V=\H^1_{0}(\Omega)$. When we consider the boundary condition $\varphi=0$ on $\partial\Omega$ in problem \eqref{PFort}, a similar construction can be done. However, we need a special care to treat the compatibility conditions in $\Omega_{-}$. Starting from the same Ansatz \eqref{1E6}, we get
\begin{equation}
\left\{
   \begin{array}{lll}
   \Delta \varphi_{0}^-& =\quad0 \quad
   &\mbox{in}\quad \Omega_{-}
\\[0.5ex]
   \partial_{\nn}\varphi_{0}^-& =\quad g \quad
   &\mbox{on}\quad \Sigma
   \end{array}
    \right.
 \label{3E40}
\end{equation}
and
\begin{equation}
\left\{
   \begin{array}{lll}
   \Delta \varphi_{0}^+ & = \quad f_{+} \quad
   &\mbox{in}\quad \Omega_{+}
\\[0.5ex]
   \varphi_{0}^+ & = \quad\varphi_{0}^-  \quad &\mbox{on}\quad \Sigma
\\[0.5ex]
   \varphi_{0}^+&=\quad 0 \quad
   &\mbox{on}\quad \partial\Omega
   \end{array}
    \right.
 \label{3E41}
\end{equation}
and for $k=1,2,...$
\begin{equation}
\left\{
   \begin{array}{lll}
   \Delta \varphi_{k}^-&=\quad\delta_{k}^1\, f^{-}\quad
   &\mbox{in}\quad \Omega_{-}
\\[0.5ex]
   \partial_{\nn}\varphi_{k}^-&=\quad -\delta_{k}^1\, g +\partial_{\nn}\varphi_{k-1}^+
   &\mbox{on}\quad \Sigma
   \end{array}
    \right.
 \label{3E42}
\end{equation}
and
\begin{equation}
\left\{
   \begin{array}{lll}
   \Delta \varphi_{k}^+ & =\quad 0 \quad
   &\mbox{in}\quad \Omega_{+}
\\[0.5ex]
   \varphi_{k}^+  &= \quad \varphi_{k}^-  \quad &\mbox{on}\quad \Sigma
\\[0.5ex]
   \varphi_{k}^+&=\quad 0 \quad
   &\mbox{on}\quad \partial\Omega.
   \end{array}
    \right.
 \label{3E43}
\end{equation}

\subsubsection{Discussion of elementary problems.}

Let $\tilde{\varphi}_{0}^-\in\H^{\frac32}(\Omega_-,\Delta)$ be the solution of the Neumann problem \eqref{3E40} under the condition $\int_{\Omega_{-}} \tilde{\varphi}_{0}^-\, \dr\xx =0$. Here we still keep a constant to be adjusted; we call it $c_{0}$. Then $\varphi_{0}\cd= \tilde{\varphi}_{0}\cd+c_{0}$ will be determined once $c_{0}$ is fixed and \eqref{3E41} will give a unique $\varphi_{0}^+\in\H^{\frac32}(\Omega_+,\Delta)$.

We consider now \eqref{3E42} for $k=1$, which is a Neumann problem with the compatibility condition
$$
   \int_{\Omega\con} f^{-}\, \dr\xx + \int_{\Sigma}(-g+\partial_{\nn}\varphi_{0}^+)\, \dr s =0
$$
and since $ \int_{\Sigma}g\, \dr s =0$, it reads
 \begin{equation}
 \label{3E50}
   \int_{\Sigma}\partial_{\nn}\varphi_{0}^+\, \dr s =- \int_{\Omega\con} f^{-}\, \dr\xx.
\end{equation}
But now we choose  $\varphi_{0}^+=\tilde\varphi_{0}^+ +c_{0}\psi$ where
\begin{equation}
\left\{
   \begin{array}{lll}
   \Delta \tilde\varphi_{0}^+ &=\quad f_{+} \quad
&\mbox{in}\quad \Omega\iso
\\[0.5ex]
   \tilde{\varphi}_{0}^+ &=\quad  \tilde{\varphi}_{0}^-  \quad &\mbox{on}\quad \Sigma
\\[0.5ex]
   \tilde{\varphi}_{0}^+ &=\quad 0 \quad
   &\mbox{on}\quad \partial\Omega.
   \end{array}
    \right.
 \label{1E45}
\end{equation}
and
\begin{equation}
\label{3E44}
\left\{
   \begin{array}{lll}
   \Delta \psi  &=\quad 0 \quad
   &\mbox{in}\quad \Omega\iso
\\[0.5ex]
   \psi &=\quad 1 \quad &\mbox{on}\quad \Sigma
\\[0.5ex]
   \psi &=\quad 0 \quad
   &\mbox{on}\quad \partial\Omega.
   \end{array}
    \right.
\end{equation}
with
\begin{equation}
   c_{0}= -\Big(\int_{\Sigma} \partial_{\nn}\tilde\varphi_{0}\is\, \dr s +
   \int_{\Omega_{-}} f^{-} \dr\xx\Big)\, \Big/\,\int_{\Sigma} \partial_{\nn}\psi\, \dr s.
\label{3E45}
\end{equation}
Clearly $\tilde\varphi_{0}\is \in \H^{\frac32}(\Omega_+,\Delta)$ is uniquely determined by \eqref{1E45}, $\psi$ by \eqref{3E44}, and $c_{0}$ by \eqref{3E45} since $\int_{\Sigma}\partial_{\nn}\psi\, \dr s\neq 0$ (principle of maximum).

Thus we have completely determined $\varphi_{0}\cd= \tilde{\varphi}_{0}\cd+c_{0} \in \H^{\frac32}(\Omega_-,\Delta)$ and $\varphi_{0}^+=\tilde\varphi_{0}^+ +c_{0}\psi \in\H^{\frac32}(\Omega_+,\Delta)$. Moreover the choice of $c_{0}$ gives the compatibility condition \eqref{3E50} of problem \eqref{3E42} for $k=1$. 

Again we take $\varphi_{1}\cd= \tilde{\varphi}_{1}\cd+c_{1} \in \H^{\frac32}(\Omega_-,\Delta)$, with $\tilde{\varphi}_{1}\cd$ uniquely determined by \eqref{3E42} with $k=1$ under the condition $\int_{\Omega_{-}} \tilde{\varphi}_{1}^-\, \dr\xx =0$. Then
$\varphi_{1}^+=\tilde\varphi_{1}^+ +c_{1}\psi \in \H^{\frac32}(\Omega_+,\Delta)$, with $\tilde\varphi_{1}^+$ uniquely determined by \eqref{3E43} with $\tilde\varphi_{1}^+=\tilde\varphi_{1}^-$ on $\Sigma$ and
\begin{equation}
c_{1}= -\int_{\Sigma} \partial_{\nn}\tilde\varphi_{1}\is\, \dr s\, \Big/\,\int_{\Sigma} \partial_{\nn}\psi\, \dr s.
\label{3E46}
\end{equation}
We can continue this iterative process to construct the sequences
 $\{\varphi_{k}\cd\}_{k\geqslant 0}\subset \H^{\frac32}(\Omega_-,\Delta)$ and $\{\varphi_{k}^+\}_{k\geqslant 0}\subset \H^{\frac32}(\Omega_+,\Delta)$.

\subsubsection{Proof of Theorem \ref{2T9} in the Dirichlet case} 
The absolute convergence of the series $\sum_{k\geqslant0}\rho^{-k}\varphi^\pm_k$ in $\H^{\frac32}(\Omega_\pm,\Delta)$ is obtained like in the Neumann case. The proof of the uniqueness of solutions to problem \eqref{1E1} in the Dirichlet case is also similar to the Neumann case.

\subsection{Complements}
\label{S3.3}

In this subsection we give some complementary results, sharper or more general than those of Theorem \ref{2T9}.

\subsubsection{Uniform estimates for $\varphi-\varphi_{0}$}  
As a consequence of the bounds \eqref{eq:3.20}, we have
\[
   \varphi-\varphi_{0}=\sum_{n=1}^\infty\, \varphi_{n} \,\rho^{-n}
\]
and we deduce the following estimate between the solution $\varphi\in V$ of the problem \eqref{1E1} and the solution $\varphi_{0}$ of the limit problem as $\rho$ tends to infinity.

\begin{thm}
\label{3T1}
Let us assume that $a_+\neq0$.
There exist a constant $\rho_{0}>0$ independent of $a_+$ such that for all $a_-\in\{z\in\C \, | \,|z|\geqslant \rho_{0}|a_+| \}$, the unique solution $\varphi\in V$ of the problem \eqref{1E1} with data $(f,g)$ satisfying \eqref{1E7}-\eqref{1E2} converges in the piecewise $\H^{3/2}$ norm to the solution $\varphi_{0}$ of the limit problem as $\rho$ tends to infinity,  with the uniform estimate
\begin{equation}
\label{ue0}
   \|\varphi^+ -\varphi^+_{0}\|_{\frac32,\Omega_+} + \|\varphi^- -\varphi^-_{0}\|_{\frac32,\Omega_-}
   \leqslant C_{\rho_{0}}  |\rho|^{-1} \big( |a_+|^{-1} \|f \|_{0,\Omega} + \|g\|_{0,\Sigma}\big)
   \end{equation}
with a constant $C_{\rho_{0}}>0$, independent of $a_+$, $a_-$, $f$, and $g$.
\end{thm}
In this context, the above result gives sharper estimates than \cite{HNM09} where we find a characterization of limit solutions and strong convergence results for similar (and more general) problems.

Likewise, an estimate of the remainder at any order is valid:
\begin{equation*}
   \|\varphi^+ -\sum_{n=0}^K\, \varphi^+_{n} \rho^{-n}\|_{\frac32,\Omega_+} 
   + \|\varphi^- -\sum_{n=0}^K\, \varphi^-_{n} \rho^{-n}\|_{\frac32,\Omega_-} \!\!
   \leqslant C_{\rho_{0}}  |\rho|^{-1-K} \big( |a_+|^{-1} \|f \|_{0,\Omega} + \|g\|_{0,\Sigma}\big)
   \end{equation*}

\subsubsection{Uniform estimates when $a_-$ and $a_+$ play symmetric roles} 
The framework of Theorem \ref{2T9} can be extended so that $a_-$ and $a_+$ play symmetric roles, and so that the contributions of the norms $\|f ^-\|_{0,\Omega_-}$ and $\|f ^+\|_{0,\Omega_+}$ are optimally taken into account in the right hand side of estimates. For this, we require the following assumptions
\begin{subequations}
\label{3H}
\begin{eqnarray}
\label{3Ha}
   \int_{\Omega}  f\, \dr\xx = \int_{\Sigma} g \, \dr s = 0 
   &&\mbox{if \ $V=V_\NN$} \\
\label{3Hb}
   \qquad\int_{\Omega_+}\!\!  f^+\, \dr\xx = 
   \int_{\Omega_-}\!\!  f^-\, \dr\xx =  \int_{\Sigma} g \, \dr s = 0 
   &&\mbox{if \ $V=V_\DD$ \ and \ $|a_-| \gg |a_+|$} \\
\label{3Hc}
   \mbox{no condition}
   &&\mbox{if \ $V=V_\DD$ \ and \ $|a_-| \ll |a_+|$}.
\end{eqnarray}
\end{subequations}

\begin{prop}
\label{3P1}
Let us assume that $a_\pm\neq0$.
There exist a constant $\rho_{0}>0$ such that for all couples $(a_-,a_+)$ such that
\[
   |a_-| \geqslant \rho_{0}|a_+| 
   \quad\mbox{or}\quad
   |a_+| \geqslant \rho_{0}|a_-|
\] 
the problem \eqref{1E1} with data $(f,g)$ satisfying \eqref{3H} has a unique solution $\varphi\in V$, which moreover is piecewise $\H^{3/2}$ and satisfies the uniform estimate
\begin{equation}
\label{ue3P1}
   \|\varphi^+\|_{\frac32,\Omega_+} + \|\varphi^-\|_{\frac32,\Omega_-}
   \leqslant C_{\rho_{0}} \bigg(
   \frac{\|f^- \|_{0,\Omega_-}} {|a_-|} + \frac{\|f^+ \|_{0,\Omega_+}} {|a_+|} + 
   \|g\|_{0,\Sigma}
   \bigg)
\end{equation}
with a constant $C_{\rho_{0}}>0$, independent of $a_+$, $a_-$, $f$, and $g$.
\end{prop}

\begin{proof}
1) \ Let us first prove estimate \eqref{ue3P1} in the Neumann case and when the modulus of $\rho:= a_-(a_+)^{-1}$ is large enough. After the change of data $(f,g)\to (a_+^{-1} f,g)$ as explained at the beginning of subsection \ref{S3.0}, proving estimate \eqref{ue3P1} reduces to show
\begin{equation}
\label{1E56b}
   \|\varphi^+_{(\rho)}\|_{\frac32,\Omega_+} + \|\varphi^-_{(\rho)}\|_{\frac32,\Omega_-}
   \leqslant C{(\rho_{0})}\, \big(
   |\rho^{-1}|\, \|f^- \|_{0,\Omega_-} +  \|f^+ \|_{0,\Omega_+} + \|g\|_{0,\Sigma}\big),
\end{equation}
instead of \eqref{1E56} (note that the new factor $|\rho^{-1}|$ in front of $\|f^- \|_{0,\Omega_-}$ is equal to $|a_+|/|a_-|$).

Estimate \eqref{1E56b} is in fact a mere consequence of estimates \eqref{1E52} and \eqref{1E55}, where we take advantage of the presence of the factor $|\rho^{-1}|$ in front of the norm $\|\varphi_{1}^- \|_{\H^{\frac32}(\Omega_-,\Delta)}$ of $\varphi^-_1$ in \eqref{1E52}.

\smallskip\noindent
2) \ Still in the Neumann case, but when $\rho:= a_+(a_-)^{-1}$ is large enough, the sequence of problems to be solved is now
\begin{equation}
\left\{
   \begin{array}{lll}
   \Delta \varphi_{0}^+ & = \quad 0 \quad
   &\mbox{in}\quad \Omega_{+}
\\[0.5ex]
   \partial_\nn\varphi_{0}^+ & = \quad g  \quad &\mbox{on}\quad \Sigma
\\[0.5ex]
   \partial_\nn\varphi_{0}^+&=\quad 0 \quad
   &\mbox{on}\quad \partial\Omega
   \end{array}
    \right.
 \label{3E41b}
\end{equation}
\begin{equation}
\left\{
   \begin{array}{lll}
   \Delta \varphi_{0}^-& =\quad f^- \quad
   &\mbox{in}\quad \Omega_{-}
\\[0.5ex]
   \varphi_{0}^-& =\quad \varphi^+_0 \quad
   &\mbox{on}\quad \Sigma
   \end{array}
    \right.
 \label{3E40b}
\end{equation}
and for $k=1,2,...$
\begin{equation}
\left\{
   \begin{array}{lll}
   \Delta \varphi_{k}^+ & =\quad \delta_{k}^1\, f^{+} \quad
   &\mbox{in}\quad \Omega_{+}
\\[0.5ex]
   \partial_\nn\varphi_{k}^+  &= \quad -\delta_{k}^1\, g +\partial_\nn\varphi_{k-1}^- \quad 
   &\mbox{on}\quad \Sigma
\\[0.5ex]
   \partial_\nn\varphi_{k}^+&=\quad 0 \quad
   &\mbox{on}\quad \partial\Omega.
   \end{array}
    \right.
 \label{3E43b}
\end{equation}
\begin{equation}
\left\{
   \begin{array}{lll}
   \Delta \varphi_{k}^-&=\quad 0 \quad
   &\mbox{in}\quad \Omega_{-}
\\[0.5ex]
   \varphi_{k}^-&=\quad \varphi_{k}^+
   &\mbox{on}\quad \Sigma
   \end{array}
    \right.
 \label{3E42b}
\end{equation}
The compatibility of the right hand sides of problems \eqref{3E41b} and \eqref{3E43b} in $\Omega_+$ can be checked by arguments similar to those used in the case when $a_-(a_+)^{-1}$ is large (\S\,\ref{S321}). The estimate can be proved similarly.

\smallskip\noindent
3) \ In the Dirichlet case, if $|a_-| >>|a_+|$, under assumption \eqref{3Hb}, we see that in \eqref{3E45}, we simply have
\begin{equation*}
   c_{0}= -\int_{\Sigma} \partial_{\nn}\tilde\varphi_{0}\is\, \dr s \; \Big/\,
   \int_{\Sigma} \partial_{\nn}\psi\, \dr s.
\end{equation*}
Thus $f^-$ does not influence $\varphi_0$, and we have estimates like in \eqref{1E52} with the factor $|\rho^{-1}|$ in front of $\|\varphi_{1}^- \|_{\H^{\frac32}(\Omega_-,\Delta)}$. We deduce estimate \eqref{ue3P1} like in the Neumann case.

\smallskip\noindent
4) \ Finally, in the Dirichlet case, if $|a_-| <<|a_+|$, none of the elementary problems is of Neumann type. Hence no compatibility condition is required and we can prove estimate \eqref{ue3P1} as previously.
\end{proof}

The compatibility conditions \eqref{1E2} (and a fortiori \eqref{3Ha}-\eqref{3Hb}) are not necessary for the solvability of problem \eqref{1E1}: For Neumann exterior boundary condition, the necessary and sufficient condition is
\[
 - \int_{\Omega}  f\, \dr\xx + (a_+-a_-)\int_{\Sigma} g \, \dr s = 0.
\]
It depends on coefficients $a_\pm$. If we want to have the compatibility of the right hand side for any value of the coefficients $a_\pm$ we can either assume \eqref{1E2} or replace the coefficient in front of the integral $\int_\Sigma$ by $1$, defining the new problem
\begin{equation}
   \forall \psi \in  V, \;
   \int_{\Omega_{+}} \!\! a_+\nabla\varphi^+ \!\cdot\, \nabla\overline{\psi}{}^+ \dr\xx
   + \int_{\Omega_{-}} \!\! a_-\nabla\varphi^- \!\cdot\, \nabla\overline{\psi}{}^- \dr\xx= 
   -\!\int_{\Omega}  f \; \overline{\psi}\,\dr\xx 
   + \!\int_{\Sigma} g\; \overline{\psi} \, \dr s. \hskip-1ex
\label{1E1b}
\end{equation}

\begin{prop}
\label{3P2}
If we assume the compatibility conditions 
\begin{equation}
 -\int_{\Omega}  f\, \dr\xx +
 \int_{\Sigma} g \, \dr s = 0\ \mbox{ if } \ V=V_\NN \quad \mbox{and} \quad 
 \mbox{nothing} \ \mbox{ if } \ V= V_\DD\,,
\label{1E2c}
\end{equation}
then problem \eqref{1E1b} is uniquely solvable if the modulus of $\rho:=a_-(a_+)^{-1}$ is large enough and its solution satisfies the uniform estimate 
\begin{equation*}
   \|\varphi^+\|_{\frac32,\Omega_+} + \|\varphi^-\|_{\frac32,\Omega_-}
   \leqslant C_{\rho_{0}}  {|a_+|}^{-1} \big(
   {\|f \|_{0,\Omega}}  +   \|g\|_{0,\Sigma}
   \big)\ .
\end{equation*}
\end{prop}

\begin{rem}
We pay the weaker assumption on the data by a weaker estimate since we have the same estimate for solutions of problem \eqref{1E1b} as for solutions of problem \eqref{1E1}: In problem \eqref{1E1} the interface datum is $(a_+-a_-)g$ (the coefficient $|a_+-a_-|$ tends to infinity as $|\rho|\to\infty$) whereas the interface datum of problem \eqref{1E1b} is $g$ alone.
\end{rem}

\begin{proof}
The construction of the terms of the series expansion is similar as in the proof of Theorem \ref{2T9}. Now we have $\varphi^-_0=0$ in the Neumann case, and $\varphi^-_0=c_0$ in the Dirichlet case.
\end{proof}

\section{Proof of uniform estimates for Maxwell solutions at high conductivity}
\label{S4}

We consider now the harmonic Maxwell system \eqref{0E2} at a fixed frequency $\omega$ satisfying Hypothesis \ref{H1}. We are going to prove the following sequence of statements:

\begin{lem}
\label{5L1}
Under Hypothesis {\em\ref{H1}}, there are constants $\sigma_{0}$ and $C_0>0$ such that if $\sigma\geqslant\sigma_0$ any solution $(\EE,\HH)\in\bL^2(\Omega)^2$ of problem \eqref{0E2} with boundary condition \eqref{ME2} and data $\jj\in\bH_0(\Div,\Omega)$ satisfies the estimate
\begin{equation}
   \|\EE\|_{0,\Omega} \leqslant C_0 \|\jj\|_{\bH(\Div,\Omega)}.
\label{5E1}
\end{equation}
A similar statement holds for boundary conditions \eqref{ME1} and data $\jj\in\bH(\Div,\Omega)$.
\end{lem}

This lemma is the key for the proof of Theorem \ref{4T0} and is going to be proved in the next subsection, using in particular our uniform estimates in the scalar case (this is the main difference with the proof of Theorem 2.1 in \cite{H-J-N07}). As a consequence of this lemma, we will obtain estimates \eqref{3E1}:

\begin{cor}
\label{5C1}
Let $\sigma>0$. Let $(\EE,\HH)\in\bL^2(\Omega)^2$ be solution of problem \eqref{0E2} with boundary condition \eqref{ME2} and data $\jj\in\bH_0(\Div,\Omega)$.
If $\EE$ satisfies estimate \eqref{5E1}, then setting 
\[
   C_1 = 1 + (1+\sqrt{\omega\mu_0})\sqrt{C_0} + (1+\omega\sqrt{\varepsilon_0\mu_0})C_0\,,
\]
there holds
\begin{equation}
   \|\EE\|_{0,\Omega} + \|\rot\EE\|_{0,\Omega} + 
   \|\Div(i\omega\varepsilon_0 - \bs) \EE\|_{0,\Omega}
   + \sqrt\sigma\, \|\EE\|_{0,\Omega_-} \leqslant C_1 \|\jj\|_{\bH(\Div,\Omega)}.
\label{5E2}
\end{equation}
A similar estimate holds for boundary conditions \eqref{ME1} and data $\jj\in\bH(\Div,\Omega)$.
\end{cor}

Finally, estimate \eqref{5E1} implies existence and uniqueness of solutions.
\begin{cor}
\label{5C2}
Let $\sigma>0$. We assume that estimate \eqref{5E1} holds for any solution $(\EE,\HH)\in\bL^2(\Omega)^2$ of problem \eqref{0E2}-\eqref{ME2} with $\jj\in\bH_0(\Div,\Omega)$. 
Then for any $\jj\in\bH_0(\Div,\Omega)$, there exists a unique solution $(\EE,\HH)\in\bL^2(\Omega)^2$ of problem \eqref{0E2}-\eqref{ME2}.
A similar result holds for boundary conditions \eqref{ME1} and data $\jj\in\bH(\Div,\Omega)$.
\end{cor}

The previous three statements clearly imply Theorem \ref{4T0}. To prepare for their proofs, we recall variational formulations in electric field for the Maxwell problem \eqref{0E2} with boundary condition \eqref{ME1} or \eqref{ME2}, cf.\ \cite{Mo03} for instance. Let
\begin{subequations}
\begin{gather}
   \bH(\rot,\Omega) = \{\uu\in\bL^2(\Omega)\; |\quad \rot\uu\in\bL^2(\Omega) \} \\
   \bH_0(\rot,\Omega) = \{\uu\in\bH(\rot,\Omega)\; |\quad 
   \uu\times\nn=0 \ \ \mbox{ on } \ \partial\Omega \}.
\end{gather}
\end{subequations}
If $(\EE,\HH)\in\bL^2(\Omega)^2$ is solution of \eqref{0E2}-\eqref{ME2}, then $\EE\in\bH(\rot,\Omega)$ satisfies for all $\EE'\in\bH(\rot,\Omega)$:
\begin{equation}
\label{5EV}
   \int_{\Omega} \big(\rot\EE \cdot \rot\overline{\EE'}  
   - \kappa^2 \EE \cdot \overline{\EE'} \big)\;\dr\xx 
   - i\nu\sigma \int_{\Omega_-} \!\! \EE \cdot \overline{\EE'}\;\dr\xx 
   = i\nu  \int_{\Omega} \jj\cdot \overline{\EE'}\,\dr\xx
\end{equation}
where we have set $\kappa=\omega\sqrt{\varepsilon_0\mu_0}$ and $\nu=\omega\mu_0$.
If boundary conditions \eqref{ME1} are considered, then $\EE\in\bH_0(\rot,\Omega)$ and \eqref{5EV} holds for any $\EE'\in\bH_0(\rot,\Omega)$.

\subsection{Proof of Lemma \ref{5L1}\,: Uniform $\bL^2$ estimate of the electric field}
Reductio ad absurdum: We assume that there is a sequence $(\EE_m,\HH_m)\in\bL^2(\Omega)^2$, $m\in\N$, of solutions of the Maxwell system \eqref{0E2}-\eqref{ME2} associated with a conductivity $\sigma_m$ and a right hand side $\jj_m\in\bH_0(\Div,\Omega)$:
\begin{subequations}
\begin{gather}
\label{5E3a}
    \rot  \EE_m - i \omega\mu_0 \HH_m = 0 \quad\mbox{in}\quad\Omega\,,\\
\label{5E3b}
    \rot \HH_m + (i\omega\varepsilon_0 - \bs_m) \EE_m = \jj_m \quad\mbox{in}\quad\Omega\,,\\
\label{5E3c}
   \HH_m\times\nn = 0 \quad\mbox{on}\quad\partial\Omega\,,
\end{gather}
\end{subequations}
satisfying the following conditions
\begin{subequations}
\begin{eqnarray}
\label{5E4a}
   &\sigma_m\to\infty\quad &\mbox{as \ $m\to\infty$,} \\
\label{5E4b}
   &\|\EE_m\|_{0,\Omega} = 1\quad &\mbox{$\forall m\in\N$,} \\
\label{5E4c}
   &\|\jj_m\|_{\bH(\Div,\Omega)}\to0\quad&\mbox{as \ $m\to\infty$.}
\end{eqnarray}
\end{subequations}
Note that the external boundary condition $\EE_m\cdot\nn=0$ on $\partial\Omega$ is but a consequence of the equation \eqref{5E3b}, the boundary condition \eqref{5E3c} and the condition $\jj\cdot\nn=0$ on $\partial\Omega$ contained in the assumption that $\jj_m$ belongs to $\bH_0(\Div,\Omega)$.

We particularize the electric variational formulation \eqref{5EV} for the sequence $\{\EE_m\}$: For all $\EE'\in\bH(\rot,\Omega)$:
\begin{equation}
\label{5EVm}
   \int_{\Omega} \big(\rot\EE_m \cdot \rot\overline{\EE'}  
   - \kappa^2 \EE_m \cdot \overline{\EE'} \big)\;\dr\xx 
   - i\nu\sigma_m \int_{\Omega_-} \!\! \EE_m \cdot \overline{\EE'}\;\dr\xx 
   = i\nu  \int_{\Omega} \jj_m\cdot \overline{\EE'}\,\dr\xx\,.
\end{equation}
Choosing $\EE'=\EE_m$ in \eqref{5EVm} and taking the real part, we obtain with the help of condition \eqref{5E4b} the following uniform bound on the curls
\begin{equation}
\label{5EC}
   \| \rot\EE_m \|_{0,\Omega} \leqslant \kappa + \sqrt{\nu \| \jj_m \|_{0,\Omega}}\,.
\end{equation}

\subsubsection{Decomposition of the electric field and bound in $\H^{\frac12}$}
We recall that we have assumed that the domain $\Omega$ is simply connected and has a smooth connected boundary.
Relying to Theorem 2.9 and Theorem 3.12 in \cite{ABDG98}, we obtain that
for all $n\in\N$ there exists a unique $\ww_m\in\H^1(\Omega)^3$ such that 
\begin{equation}
\label{5E5}
   \rot\ww_m=\rot\EE_m,\quad  \Div \ww_m=0 \ \mbox{ in } \ \Omega ,\quad \mbox{and} \quad \ww_m\cdot\nn=0\ \mbox{ on }\ \partial\Omega \,.
\end{equation}
Moreover, we have the estimate
\begin{equation}
\label{5E6}
   \|\ww_m\|_{1,\Omega} \leqslant C \| \rot\EE_m \|_{0,\Omega}\,,
\end{equation}
where $C$ is independent of $m$. As a consequence of the equality $\rot\ww_m=\rot\EE_m$ and the simple connectedness of $\Omega$, we obtain that there exists $\varphi_m\in\H^1(\Omega)$ such that
\begin{equation}
\label{5E7}
   \EE_m = \ww_m + \nabla \varphi_m\,.
\end{equation}
We write equation \eqref{5E3b} as
\[
    \rot \HH_m + (i\omega\varepsilon_0 - \bs_m) (\ww_m + \nabla \varphi_m) = \jj_m .
\]
Let $\psi\in\H^1(\Omega)$ be a test function. Multiplying the above equality by $\nabla \overline\psi$ and integrating over $\Omega$, we obtain, using that $\Div\ww_m=0$:
\begin{equation}
\label{5E8}
   \int_{\Omega}  (i\omega\varepsilon_0 - \bs_m)\;\nabla\varphi_m \cdot 
   \nabla\overline{\psi}\;  \dr\xx = 
   - \int_{\Omega}  \Div\jj_m \; \overline\psi\;\dr\xx 
   - \sigma_m \int_{\Sigma}  \ww_m\cdot\nn\on{\Sigma} \; \overline \psi \;\dr s\,.
\end{equation}
Note that the boundary values $\rot\HH_m\cdot\nn = 0$ and $\jj_m\cdot\nn = 0$ on $\partial\Omega$ have been used here.

Thus $\varphi_m$ is solution of the Neumann problem defined by the variational equation \eqref{5E8}. Since 
\[
   \Div\jj_m\in\L^2(\Omega) \quad\mbox{and}\quad
   \int_\Omega  \Div\jj_m \,\dr\xx = \int_{\partial\Omega} \jj_m\cdot\nn  \;\dr s = 0,
\]
and 
\[
   \ww_m\cdot\nn\in\L^2(\Sigma) \quad\mbox{and}\quad
   \int_{\Sigma}  \ww_m\cdot\nn\on{\Sigma} \;\dr s = \int_{\Omega_-} \Div\ww_m \,\dr\xx = 0,
\]
the Neumann problem defined by \eqref{5E8} satisfies the assumptions of Theorem \ref{2T9} with $a_-=i\omega\varepsilon_0 - \sigma_m$ and $a_+ = i\omega\varepsilon_0$. Therefore we have the following uniform estimate for $\sigma_m$ large enough (i.e.\ for $m$ large enough, cf.\ \eqref{5E4a})
\begin{equation*}
   \|\varphi^+_m\|_{\frac32,\Omega_{+}} + \|\varphi^-_m\|_{\frac32,\Omega_{-}}\leqslant 
   C_{0} \big( \|\Div\jj_m \|_{0,\Omega}+ \|\ww_m \cdot \nn\|_{0,\Sigma} \big).
\end{equation*}
Since $\|\ww_m \cdot \nn\|_{0,\Sigma}$ is bounded by $\|\ww_m\|_{1,\Omega}$, the above inequality implies 
\begin{equation}
\label{5E9b}
   \|\varphi^+_m\|_{\frac32,\Omega_{+}} + \|\varphi^-_m\|_{\frac32,\Omega_{-}}\leqslant 
   C_{0} \big( \|\Div\jj_m \|_{0,\Omega}+ \|\ww_m\|_{1,\Omega} \big).
\end{equation}
Finally \eqref{5E4c}, \eqref{5EC}, \eqref{5E6} and \eqref{5E9b} implies that
\begin{equation}
\label{5E11}
   \|\varphi^+_m\|_{\frac32,\Omega_{+}} + \|\varphi^-_m\|_{\frac32,\Omega_{-}}
   + \|\ww_m\|_{1,\Omega} \leqslant B
\end{equation}
for a constant $B>0$ independent of $m$. With \eqref{5E7}, \eqref{5E11} gives that the sequence $\{\EE_m\}$ is bounded in $\H^{\frac12}$ on $\Omega_-$ and $\Omega_+$:
\begin{equation*}
   \|\EE^+_m\|_{\frac12,\Omega_{+}} + \|\EE^-_m\|_{\frac12,\Omega_{-}}
   \leqslant B .
\end{equation*}
Combining the above bound with \eqref{5EC}, we obtain the uniform bound
\begin{equation}
\label{5E12}
   \|\EE^+_m\|_{\frac12,\Omega_{+}} + \|\EE^-_m\|_{\frac12,\Omega_{-}} +
   \| \rot\EE_m \|_{0,\Omega} 
   \leqslant C .
\end{equation}

\subsubsection{Limit of the sequence and conclusion}
The domains $\Omega_\pm$ being bounded, the embedding of $\H^{\frac12}(\Omega_\pm)$ in $\L^2(\Omega_\pm)$ is compact. Hence as a consequence of \eqref{5E12}, we can extract a subsequence  of $\{\EE_m\}$ (still denoted by $\{\EE_m\}$) which is converging in $\bL^2(\Omega)$. By the Banach-Alaoglu theorem, we can assume that the sequence $\rot\EE_m$ is weakly converging in $\bL^2(\Omega)$: We deduce that there is $\EE\in \bL^2(\Omega)$ such that
\begin{equation}
\left\{
   \begin{array}{lll}
   \label{4E6}
   \rot\EE_m \rightharpoonup \rot\EE   \quad &\mbox{in} \quad \bL^2(\Omega)
\\
   \EE_m \rightarrow \EE   \quad &\mbox{in} \quad \bL^2(\Omega).
   \end{array}
\right.
\end{equation}
A consequence of the strong convergence in $\bL^2(\Omega)$ and \eqref{5E4b} is that $\|\EE\|_{0,\Omega}=1$.
Using Hypothesis \eqref{H1}, we are going  to prove that $\EE=0$, which will contradict $\|\EE\|_{0,\Omega}=1$, and finally prove estimate \eqref{5E1}.

Taking imaginary parts in \eqref{5EVm} when $\EE_m$ is the test-function, then letting $m\rightarrow + \infty$ and using \eqref{5E4c} we get  $\|\EE\|_{0,\Omega_{-}} =0$. Hence,
\begin{equation}
   \EE=0 \quad \mbox{in} \quad \Omega_{-}\,.
\label{4E9}
\end{equation}

Let us introduce the space
\[
   \bH_0(\rot,\Omega_{+},\Sigma) := \{\uu\in\bH(\rot,\Omega_{+})\;|\quad 
   \uu\times\nn=0 \ \mbox{ on }\ \Sigma\}.
\]
In particular, \eqref{4E9} implies that $\EE^+:=\EE\on{\Omega_+}$ belongs to $\bH_0(\rot,\Omega_{+},\Sigma)$.

Let $\Phi \in \bH_0(\rot,\Omega_{+},\Sigma)$. Then the extension $\Phi_0$ of $\Phi$ by $0$ on $\Omega_-$ defines an element of $\bH(\rot,\Omega)$. We can use $\Phi_0$ as test function in \eqref{5EVm} and we obtain
\begin{equation*}
   \int_{\Omega_+} \big(\rot\EE_m \cdot \rot\overline{\Phi}  
   - \kappa^2 \EE_m \cdot \overline{\Phi} \big)\;\dr\xx 
   = i\nu  \int_{\Omega_+} \jj_m\cdot \overline{\Phi}\,\dr\xx\,.
\end{equation*}
According to \eqref{4E6} and \eqref{5E4c}, taking limits as $m\rightarrow +\infty$, we deduce from the previous equalities
\begin{equation}
\label{5E20}
   \int_{\Omega_+} \big(\rot\EE^+ \cdot \rot\overline{\Phi}  
   - \kappa^2 \EE^+ \cdot \overline{\Phi} \big)\;\dr\xx = 0\,,
\end{equation}
i.e., $\EE^+ \in \bH_0(\rot,\Omega_{+},\Sigma)$ satisfies \eqref{5E20} for all $\Phi \in \bH_0(\rot,\Omega_{+},\Sigma)$. 
Integrating by parts we find (with $\Phi_\TT$ the tangential part of $\Phi$ on $\partial\Omega$)
\begin{equation*}
   (\rot\EE^+, \rot\Phi)_{0,\Omega\iso} = 
   (\rot\rot\EE^+, \Phi)_{0,\Omega\iso} - 
   (\rot\EE^+\times\nn, \Phi_{\TT})_{|\partial\Omega}.
\end{equation*}
Thus we have
\begin{equation}
\left\{
   \begin{array}{lll}
   \rot\rot\EE^+ - \kappa^2\EE^+ &=\quad 0 \ \quad & \mbox{in}\quad \Omega\iso
\\[0.5ex]
   \EE^+\times\nn &=\quad 0 \quad  &\mbox{on}\quad \Sigma
\\[0.5ex]
   \rot\EE^+\times\nn &=\quad 0  \quad &\mbox{on}\quad \partial\Omega.
   \end{array}
    \right.
 \label{4E7}
\end{equation}
Setting $\HH^+:=(i\omega\mu_0)^{-1}\rot\EE^+$, we obtain that $\rot\HH^+=-i\omega\varepsilon_0\EE^+$ and we deduce the remaining boundary conditions
\[
   \HH^+\cdot\nn=0 \ \mbox{ on }\ \Sigma
   \quad\mbox{and}\quad
   \EE^+\cdot\nn=0 \ \mbox{ on }\ \partial\Omega
\]
from the previous relations. Hence $(\EE^+,\HH^+)\in\bL^2(\Omega_+)^2$ is solution of problem \eqref{PH1}. By Hypothesis \ref{H1}, we deduce
\begin{equation*}
  \EE^+=0 \quad \mbox{in} \quad \Omega\iso.
\end{equation*}
Hence, with \eqref{4E9}, we have 
$\EE=0$ in $\Omega$,
which contradicts $\|\EE\|_{0,\Omega}=1$ and ends the proof of Lemma \ref{5L1}.

\subsection{Proof of Corollary \ref{5C1}.}
\label{S4.4}
Let $(\EE,\HH)\in\bL^2(\Omega)^2$ be a solution of the Maxwell problem \eqref{0E2} with boundary condition \eqref{ME2} and data $\jj\in\bH_0(\Div,\Omega)$. We assume that 
\begin{equation}
   \|\EE\|_{0,\Omega} \leqslant C_0 \|\jj\|_{\bH(\Div,\Omega)}.
\label{5E21}
\end{equation}
Then $\EE\in\bH(\rot,\Omega)$ is solution of the variational problem \eqref{5EV}. Taking as test function $\EE$ itself, we obtain the identity
\begin{equation}
\label{4E11}
   \int_{\Omega} \big(\rot\EE \cdot \rot\overline{\EE}  
   - \kappa^2 \EE \cdot \overline{\EE} \big)\;\dr\xx 
   - i\nu\sigma \int_{\Omega_-} \!\! \EE \cdot \overline{\EE}\;\dr\xx 
   = i\nu  \int_{\Omega} \jj\cdot \overline{\EE}\,\dr\xx\,.
\end{equation}
Taking the real part of \eqref{4E11}, we obtain
\begin{equation*}
   \|\rot\EE\|^2_{0,\Omega} = \kappa^2\|\EE\|^2_{0,\Omega} - 
   \nu\,\Im(\jj,\EE)_{0,\Omega}
\end{equation*}
hence, using inequality \eqref{5E21} and Cauchy-Schwarz inequality,
\begin{equation}
   \|\rot\EE\|_{0,\Omega}\leqslant 
   \big(\kappa C_0 + \sqrt{\nu C_0}\,\big) \| \jj\|_{\HH(\Div,\Omega)}.
\label{4E12}
\end{equation}
Then, taking the imaginary part of \eqref{4E11},
\begin{equation*}
   \sigma\|\EE\|^2_{0,\Omega_-} = -\Re(\jj,\EE)_{0,\Omega},
\end{equation*}
hence,
\begin{equation}
\label{4E13}
   \sqrt\sigma\, \|\EE\|_{0,\Omega_-} \leqslant \sqrt{C_0} \,\| \jj\|_{\HH(\Div,\Omega)}.
\end{equation}
Taking the divergence of equation $\rot\HH + (i\omega\varepsilon_0 - \bs) \EE = \jj$, we immediately obtain
\begin{equation}
\label{4E14}
   \|\Div(i\omega\varepsilon_0 - \bs) \EE\|_{0,\Omega} = 
   \|\Div\jj\|_{0,\Omega}.
\end{equation}
Formulas \eqref{4E12} to \eqref{4E14} yield Corollary \ref{5C1}.

\subsection{Proof of Corollary \ref{5C2}.}
Let $\sigma$ and $\omega$ (i.e., $\kappa$) be fixed.
Let us introduce the piecewise constant function $\bal$ on $\Omega$
\begin{equation}
\label{5E31}
   \bal = \mathds{1} +  \frac{i}{\omega\varepsilon_0} \bs \ .
\end{equation}
With this notation, the sesquilinear form in the left hand side of \eqref{5EV} becomes
\begin{equation}
\label{5E32}
   \int_{\Omega} \big(\rot\EE \cdot \rot\overline{\EE'}  
   - \kappa^2 \bal\, \EE \cdot \overline{\EE'} \big)\;\dr\xx .
\end{equation} 
The proof of Corollary \ref{5C2} relies on a classical regularization procedure:
We consider the functional space
$$
   \XX_{\TT}(\bal) = \{\EE \in \HH(\rot,\Omega) | \quad 
   \Div(\bal\EE) \in \L^2(\Omega),\quad \EE\cdot\nn = 0 \; \mbox{on}\; \partial\Omega\}.
$$
Let $s>0$ be a real number, which will be chosen later. Let us introduce the sesquilinear forms $A_s$ and $B$\,: $\XX_{\TT}(\bal)\times\XX_{\TT}(\bal)\rightarrow\C$
\begin{subequations}
\begin{gather}
   A_s(\EE,\EE') =
   \int_{\Omega} \big(\rot\EE \cdot \rot\overline{\EE'} 
   + s\, \Div\bal\EE\, \Div\overline{\bal\EE'} \big) \,\dr\xx \\
   B(\EE,\EE') =
   \int_{\Omega}   
   \bal \EE \cdot \overline{\EE^\prime}\,\dr\xx\,.
\end{gather}
\end{subequations}
With a right hand side $\jj\in\bH_0(\Div,\Omega)$, we associate a new right hand side $\jj_s$ depending on the parameter $s$ defined as an element of $\XX_\TT(\bal)'$ by
\begin{equation}
   \jj_s(\EE') = \int_{\Omega} \big( \jj \cdot\overline{\EE'} 
   -\frac{s}{\kappa^2} \Div\jj\, \Div\overline{\bal\EE'}\big) \,\dr\xx \quad  
   \forall\EE'\in\XX_\TT(\bal).
\end{equation}
The regularized variational formulation is: Find $\EE \in  \XX_{\TT}(\bal)$ such that
\begin{equation}
\label{5E35}
   \forall\,\EE' \in  \XX_{\TT}(\bal),\quad
   A_s(\EE,\EE') - \kappa^2 B(\EE,\EE') = i\nu \jj_s(\EE').
\end{equation}
As a consequence of \cite[Th.\,7.2]{CoDa00}, we obtain that if
\begin{equation}
\label{5E36}
   \frac{\kappa^2}{s} \ \ \mbox{is not an eigenvalue of the Neumann problem for the operator}
   \ \ \Div\bal\nabla,
\end{equation}
then any solution $(\EE,\HH)\in\bL^2(\Omega)^2$ of problem \eqref{0E2}-\eqref{ME2} with $\jj\in\bH_0(\Div,\Omega)$ provides a solution of problem \eqref{5E35}, and conversely, any solution $\EE$ of \eqref{5E35} provides a solution of \eqref{0E2}-\eqref{ME2} by setting $\HH=(i\omega\mu)^{-1}\rot\EE$.

Thus, we choose $s$ so that \eqref{5E36} holds.

Since the form $A_s$ is coercive on $\XX_\TT(\bal)$ and the embedding of $\bL^2(\Omega)$ in $\XX_\TT(\bal)$ is compact, we obtain that the Fredholm alternative is valid: If the kernel of the adjoint problem to \eqref{5E35}
\begin{equation}
\label{5E35k}
   \mbox{Find }\ \EE' \in  \XX_{\TT}(\bal),\quad 
   \forall\,\EE \in  \XX_{\TT}(\bal),\quad
   A_s(\EE,\EE') - \kappa^2 B(\EE,\EE') = 0,
\end{equation}
is reduced to $\{0\}$, then problem \eqref{5E35} is solvable.

We see that the assumption of Corollary \ref{5C2} implies that \eqref{5E35k} has only the zero solution, and that the same holds for the direct problem
\begin{equation}
   \mbox{Find }\ \EE \in  \XX_{\TT}(\bal),\quad 
   \forall\,\EE' \in  \XX_{\TT}(\bal),\quad
   A_s(\EE,\EE') - \kappa^2 B(\EE,\EE') = 0,
\end{equation}
of course.

All this implies the unique solvability of problem \eqref{0E2}-\eqref{ME2} with $\jj\in\bH_0(\Div,\Omega)$.

\section{Application: Convergence of asymptotic expansion at high conductivity}
\label{S5}
In the Maxwell case, see equations \eqref{0E2}, let us introduce the parameter
\begin{equation}
\label{6E1}
   \delta = \sqrt{\frac{\omega\varepsilon_0}{\sigma}}\,.
\end{equation}
Thus, when $\sigma\to\infty$, $\delta$ tends to $0$. Note that the function $\bal$ defined in \eqref{5E31} can be written
\begin{equation}
\label{6E2}
   \bal = \mathbf{1}_{{\Omega_+}}  +  
   \Big( 1+\frac{i}{\delta^2} \Big) \mathbf{1}_{\Omega_-}.
\end{equation}
Several works are devoted to the interesting question of an asymptotic expansion as $\delta\to0$ of solutions of the Maxwell system \eqref{0E2} with complementing boundary conditions on $\partial\Omega$ {\em when the interface $\Sigma$ is smooth}: See \cite{St83,McCamySte84,McCamySte85} for plane interface and eddy current approximation, \cite{H-J-N07} for impedance boundary conditions and \cite{Peron09} for perfectly insulating or perfectly conducting boundary conditions. 

\subsection{Assumptions}
\label{5S1}
We assume that $\Sigma$ is a smooth surface, and we follow the approach of \cite{Peron09}. In order to fix ideas, we take perfectly insulating boundary condition \eqref{ME2} and assume Hypothesis \ref{H1} for this condition. By Theorem \ref{4T0} there exists $\sigma_0$ such that the conclusions of the theorem hold. From now on we assume that
\begin{equation}
\label{6E3}
   \sigma\geqslant\sigma_0,\quad\mbox{i.e.}\quad \delta\leqslant\delta_0\ \ 
   \mbox{with}\ \ 
   \delta_0 = \sqrt{\frac{\omega\varepsilon_0}{\sigma_0}}\,.
\end{equation}
Let $\jj\in\bH_0(\Div,\Omega)$ such that $\jj=0$ in $\Omega_{-}$. Then for all $\delta\leqslant\delta_0$, there exists a unique solution to problem \eqref{0E2}-\eqref{ME2}, which we denote by $(\EE_{(\delta)},\HH_{(\delta)})$. Then it is possible to construct series expansions in powers of $\delta$ for the electric field $\EE^+_{(\delta)}$ in the dielectric part $\Omega_+$ and $\EE^-_{(\delta)}$ in the conducting part $\Omega_-$:
\begin{subequations}
\label{6E4}
\begin{gather}
\label{6E4a}
   \EE^+_{(\delta)}(\xx) \approx \sum_{j\geqslant0} \delta^j\EE^+_j(\xx) \\
\label{6E4b}
   \EE^-_{(\delta)}(\xx) \approx \sum_{j\geqslant0} \delta^j\EE^-_j(\xx;\delta) 
   \quad\mbox{with}\quad   
   \EE^-_j(\xx;\delta) = \chi(y_3) \,\WW_j(y_\beta,\frac{y_3}{\delta})\,.
\end{gather}
In \eqref{6E4b}, $\yy=(y_\beta,y_3)$ are ``normal coordinates'' to the surface $\Sigma$ in a tubular neighborhood $\cU_-$ of $\Sigma$ in the conductor part $\Omega_-$. In particular, $y_3$ represents the distance to $\Sigma$. The function $\yy\mapsto\chi(y_3)$ is a smooth cut-off with support in $\overline\cU_-$ and equal to $1$ in a smaller tubular neighborhood of $\Sigma$. The functions $\WW_j$ are {\em profiles} defined on $\Sigma\times\R_+$. Moreover, for any $j\in\N$
\begin{equation}
\label{6E4c}
   \EE^+_j\in\bH(\rot,\Omega_+) \quad\mbox{and}\quad
   \WW_j\in\bH(\rot,\Sigma\times\R_+).
\end{equation}
\end{subequations}
There hold a similar series expansions in powers of $\delta$ for the magnetic field $\HH_{(\delta)}$.

The validation of the asymptotic expansion \eqref{6E4} consist in proving estimates for remainders $\RR_{m;\,\delta}$ defined as
\begin{equation}
\label{6E5}
   \RR_{m;\,\delta} = \EE_{(\delta)} - \sum_{j=0}^m \delta^j\EE_j \quad\mbox{in}\quad\Omega\,.
\end{equation}
This is done by an evaluation of the right hand side when the Maxwell operator is applied to $\RR_{m;\,\delta}$. By construction \cite[Proposition 7.4]{Peron09}, we obtain
\begin{equation}
\label{6E6}
\left\{
   \begin{array}{lll}
   \rot\rot\RR^+_{m;\,\delta} - \kappa^2\alpha_+\RR^+_{m;\,\delta} &=\quad 0 
   \quad & \mbox{in}\quad \Omega_+
   \\[0.5ex]
   \rot\rot\RR^-_{m;\,\delta} - \kappa^2\alpha_-\RR^-_{m;\,\delta} &=\quad \jj^-_{m;\,\delta} 
   \quad & \mbox{in}\quad \Omega_-
   \\[0.5ex]
   \big[\RR_{m;\,\delta}\times\nn\big]_\Sigma &=\quad 0 \quad  &\mbox{on}\quad \Sigma
   \\[0.5ex]
   \big[\rot\RR_{m;\,\delta}\times\nn\big]_\Sigma &=\quad \g_{m;\,\delta}  
   \quad  &\mbox{on}\quad \Sigma
   \\[0.5ex]
   \rot\RR^+_{m;\,\delta}\times\nn &=\quad 0 \quad &\mbox{on}\quad \partial\Omega\,.
   \end{array}
    \right.
\end{equation}
Here, according to \eqref{6E2}, $\alpha_+=1$ and $\alpha_-=1+i/\delta^2$, and $[\EE\times\nn]_\Sigma$ denotes the jump of $\EE\times\nn$ across $\Sigma$. The right hand sides (residues) $\jj^-_{m;\,\delta}$ and $\g_{m;\,\delta}$ are, roughly, of the order $\delta^m$.

\subsection{Convergence result}
The main result of this section is the following.

\begin{thm}
\label{6T1}
Under Hypothesis {\rm\ref{H1}} in the framework above (section {\rm\ref{5S1}}), we assume that we have for all $m\in\N$ the following estimates for the residues $\jj^-_{m;\,\delta}$ and $\g_{m;\,\delta}$ in \eqref{6E6}
\begin{equation}
\label{6E7}
   \|\jj^-_{m;\,\delta}\|_{2,\Omega_-} + 
   \|\g_{m;\,\delta}\|_{\frac12,\Sigma} + 
   \|\rot_\Sigma\g_{m;\,\delta}\|_{\frac32,\Sigma} \leqslant
   C_m\delta^{m-m_0},
\end{equation}
where $C_m>0$ is independent of $\delta$, and $m_0\in\N$ independent of $m$ and $\delta$. Then for all $m\in\N$ and $\delta\in(0,\delta_0]$, the remainders $\RR_{m;\,\delta}$ \eqref{6E5} satisfy the optimal estimates
\begin{equation}
\label{6E8}
   \|\RR^+_{m;\,\delta}\|_{0,\Omega_+} + \|\rot\RR^+_{m;\,\delta}\|_{0,\Omega_+} + 
   \delta^{-\frac12} \|\RR^-_{m;\,\delta}\|_{0,\Omega_-} + 
   \delta^{\frac12} \|\rot\RR^-_{m;\,\delta}\|_{0,\Omega_-} \leqslant C'_m\delta^{m+1}.
   \!\!\!\!
\end{equation}
\end{thm}

\begin{proof}
Here we denote by $C_m$ various constants which may depend on $m$ but not on $\delta$.

\noindent
{\sc Step 1.} We cannot use Theorem \ref{4T0} directly because $\rot\rot\RR_{m;\,\delta} - \kappa^2\bal\RR_{m;\,\delta}$ does not define an element of $\bH(\Div,\Omega)$. We are going to introduce two correctors $\CC_{m;\,\delta}$ and $\D_{m;\,\delta}$ satisfying suitable estimates and so that
\begin{equation}
\label{6E6C}
\left\{
   \begin{array}{lll}
   \big[(\RR_{m;\,\delta} - \CC_{m;\,\delta})\times\nn\big]_\Sigma 
   &=\quad 0 \quad  &\mbox{on}\quad \Sigma
   \\[0.5ex]
   \big[\rot(\RR_{m;\,\delta}-\CC_{m;\,\delta})\times\nn\big]_\Sigma 
   &=\quad 0   \quad  &\mbox{on}\quad \Sigma
    \\[0.5ex]
   \rot(\RR_{m;\,\delta}-\CC_{m;\,\delta})\times\nn  
   &=\quad 0   \quad  &\mbox{on}\quad \partial\Omega
   \end{array}
    \right.
\end{equation}
and
\begin{equation}
\label{6E6D}
\left\{
   \begin{array}{lll}
   \big[\bal(\RR_{m;\,\delta} - \CC_{m;\,\delta} - \D_{m;\,\delta})\cdot\nn\big]_\Sigma 
   &=\quad 0 \quad  &\mbox{on}\quad \Sigma
   \\[0.5ex]
   \big[(\RR_{m;\,\delta} - \CC_{m;\,\delta} - \D_{m;\,\delta})\times\nn\big]_\Sigma 
   &=\quad 0 \quad  &\mbox{on}\quad \Sigma
   \\[0.5ex]
   \big[\rot(\RR_{m;\,\delta}-\CC_{m;\,\delta} - \D_{m;\,\delta})\times\nn\big]_\Sigma 
   &=\quad 0   \quad  &\mbox{on}\quad \Sigma
     \\[0.5ex]
   \rot(\RR_{m;\,\delta}-\CC_{m;\,\delta}- \D_{m;\,\delta})\times\nn  
   &=\quad 0   \quad  &\mbox{on}\quad \partial\Omega \ .
   \end{array}
    \right.
\end{equation}

\noindent
{\sc Step 1a.} Construction of $\CC_{m;\,\delta}$: We take $\CC_{m;\,\delta}=0$ in $\Omega_-$ and use a trace lifting to define $\CC_{m;\,\delta}$ in $\Omega_+$. It suffices that
\begin{equation}
\label{6E7C}
\left\{
   \begin{array}{lll}
   \CC^+_{m;\,\delta} \times\nn 
   &=\quad 0 \quad  &\mbox{on}\quad \Sigma
   \\[0.5ex]
   \rot\CC^+_{m;\,\delta} \times\nn 
   &=\quad \g_{m;\,\delta}   \quad  &\mbox{on}\quad \Sigma
     \\[0.5ex]
   \rot\CC^+_{m;\,\delta} \times\nn 
   &=\quad  0   \quad  &\mbox{on}\quad \partial\Omega \ .
   \end{array}
    \right.
\end{equation}
Denoting by $C_\beta$ and $C_3$ the tangential and normal components of $\CC^+_{m;\,\delta}$ associated with a system of normal coordinates $\yy=(y_\beta,y_3)$, and by $g_\beta$ the components of $\g_{m;\,\delta}$ the above system becomes (cf.\ \cite[Proposition 3.26]{Peron09})
\begin{equation}
\label{6E8C}
\left\{
   \begin{array}{lll}
   C_\beta 
   &=\quad 0 \quad  &\mbox{on}\quad \Sigma
   \\[0.5ex]
   \partial_3 C_\beta - \partial_\beta C_3 
   &=\quad g_\beta   \quad  &\mbox{on}\quad \Sigma
    \\[0.5ex]
   \partial_3 C_\beta - \partial_\beta C_3 
   &=\quad 0   \quad  &\mbox{on}\quad \partial\Omega \ .
   \end{array}
    \right.
\end{equation}
It can be solved in $\H^2(\Omega_+)$ choosing $C_3=0$ and a standard lifting of the first two traces on $\Sigma$ and $\partial\Omega$ with the estimate
\begin{equation}
\label{6E9C}
   \|\CC^+_{m;\,\delta}\|_{2,\Omega_+} \leqslant C\|\g_{m;\,\delta}\|_{\frac12,\Sigma} \ .
\end{equation}

\noindent
{\sc Step 1b.} Construction of $\D_{m;\,\delta}$: Let us denote $\RR_{m;\,\delta}-\CC_{m;\,\delta}$ by $\SS$ for short. Again, we take $\D^-_{m;\,\delta}=0$ and use a trace lifting to define $\D^+_{m;\,\delta}$. It suffices that
\begin{equation}
\label{6E7D}
\left\{
   \begin{array}{lll}
   \D^+_{m;\,\delta} \cdot \nn 
   &=\quad \big[\bal\SS\cdot\nn\big]_\Sigma 
   &\mbox{on}\quad \Sigma
   \\[0.5ex]
   \D^+_{m;\,\delta} \times\nn 
   &=\quad 0 \quad  &\mbox{on}\quad \Sigma
   \\[0.5ex]
   \rot\D^+_{m;\,\delta} \times\nn 
   &=\quad 0   \quad  &\mbox{on}\quad \Sigma\cup\partial\Omega
   \end{array}
    \right.
\end{equation}
In normal coordinates and associated components, these conditions become, compare with \eqref{6E8C}
\begin{equation}
\label{6E8D}
\left\{
   \begin{array}{lll}
   D_3 
   &=\quad \big[\bal\SS\cdot\nn\big]_\Sigma \quad  &\mbox{on}\quad \Sigma
   \\[0.5ex]
   D_\beta 
   &=\quad 0 \quad  &\mbox{on}\quad \Sigma
   \\[0.5ex]
   \partial_3 D_\beta - \partial_\beta D_3 
   &=\quad 0   \quad  &\mbox{on}\quad \Sigma\cup\partial\Omega,
   \end{array}
    \right.
\end{equation}
which can be solved in $\H^2(\Omega_+)$ (first determine $D_3$, then $D_\beta$) with the estimate
\begin{equation}
\label{6E9D}
   \|\D^+_{m;\,\delta}\|_{2,\Omega_+} \leqslant 
   C\|\big[\bal\SS\cdot\nn\big]_\Sigma\|_{\frac32,\Sigma} \ .
\end{equation}
Since $\big[\rot\SS\times\nn\big]_\Sigma=0$, we find that, by construction
\begin{equation*}
\begin{split}
   -\kappa^2 \big[\bal\SS\cdot\nn\big]_\Sigma & =
   \big[(\rot\rot\SS - \kappa^2\bal\SS)\cdot\nn\big]_\Sigma \\ &=
  - \big(\rot\rot\CC^+_{m;\,\delta} - \kappa^2\CC^+_{m;\,\delta}\big)\on\Sigma\cdot\nn 
   \ - \ \jj^-_{m;\,\delta}\on\Sigma\cdot\nn \\ &=
   \rot_\Sigma \g_{m;\,\delta} \  - \ \jj^-_{m;\,\delta}\on\Sigma\cdot\nn \ .
\end{split}
\end{equation*}
Hence
\begin{equation}
\label{6E10}
   \|\big[\bal\SS\cdot\nn\big]_\Sigma\|_{\frac32,\Sigma} \leqslant
   \|\rot_\Sigma\g_{m;\,\delta}\|_{\frac32,\Sigma} +
   \|\jj^-_{m;\,\delta}\|_{2,\Omega_-} \,.
\end{equation}
We deduce from assumption \eqref{6E7}, and \eqref{6E9C}, \eqref{6E9D},  \eqref{6E10}
\begin{equation*}
   \|\CC^+_{m;\,\delta}\|_{2,\Omega_+} + 
   \|\D^+_{m;\,\delta}\|_{2,\Omega_+} \leqslant 
   C_m\delta^{m-m_0} .
\end{equation*}
Since by construction $\CC^-_{m;\,\delta}=\D^-_{m;\,\delta}=0$  and $\big[\CC_{m;\,\delta}\times\nn\big]_\Sigma=\big[\D_{m;\,\delta}\times\nn\big]_\Sigma=0$, the above estimate implies
\begin{equation}
\label{6E11}
   \|\CC_{m;\,\delta}\|_{0,\Omega} +\|\rot\CC_{m;\,\delta}\|_{0,\Omega} + 
   \|\D_{m;\,\delta}\|_{0,\Omega} + \|\rot\D_{m;\,\delta}\|_{0,\Omega} \leqslant 
   C_m\delta^{m-m_0}\,.
\end{equation}
We set
\begin{subequations}
\begin{equation}
\label{6E12a}
   \widetilde\RR_{m;\,\delta} := \RR_{m;\,\delta} - \CC_{m;\,\delta} - \D_{m;\,\delta}
\end{equation}
and
\begin{equation}
\label{6E12b}
  \jjw_{m;\,\delta} := 
   \rot\rot\widetilde\RR_{m;\,\delta} - \kappa^2\bal\widetilde\RR_{m;\,\delta} \ .
\end{equation}
\end{subequations}
Hence by construction, $\jjw_{m;\,\delta}\in\bH(\Div,\Omega)$ with the estimates
\begin{equation}
\label{6E13}
   \|\jjw_{m;\,\delta}\|_{\bH(\Div,\Omega)}  \leqslant 
   C_m\delta^{m-m_0} .
\end{equation}

\noindent{\sc Step 2.}
We can apply Theorem \ref{4T0} to the couple $(\EE,\HH) = (\widetilde\RR_{m;\,\delta}, (i\omega\mu_0)^{-1}\rot\widetilde\RR_{m;\,\delta})$ and, thanks to \eqref{6E12b}, obtain
\begin{equation*}
   \|\widetilde\RR_{m;\,\delta}\|_{0,\Omega} + \|\rot\widetilde\RR_{m;\,\delta}\|_{0,\Omega}
   \leqslant C \|\jjw_{m;\,\delta}\|_{\bH(\Div,\Omega)}\,.
\end{equation*}
Combined with \eqref{6E13}, this gives
\begin{equation}
\label{6E14a}
   \|\widetilde\RR_{m;\,\delta}\|_{0,\Omega} + \|\rot\widetilde\RR_{m;\,\delta}\|_{0,\Omega}
   \leqslant  C_m\delta^{m-m_0}.
\end{equation}
Together with  \eqref{6E11} and \eqref{6E12a}, this estimate gives finally
\begin{equation}
\label{6E14}
   \|\RR_{m;\,\delta}\|_{0,\Omega} + \|\rot\RR_{m;\,\delta}\|_{0,\Omega}
   \leqslant  C_m\delta^{m-m_0}.
\end{equation}

\noindent{\sc Step 3.}
In order to deduce the optimal estimate \eqref{6E8} for $\RR_{m;\,\delta}$, we use \eqref{6E14} for $m+1+m_0$, which yields
\begin{equation}
\label{6E15}
   \|\RR_{m+1+m_0;\,\delta}\|_{0,\Omega} + \|\rot\RR_{m+1+m_0;\,\delta}\|_{0,\Omega}
   \leqslant  C_m\delta^{m+1}.
\end{equation}
But we have the formula
\begin{equation}
\label{6E16}
   \RR_{m;\,\delta} = \sum_{j=m+1}^{m+1+m_0} \delta^j \EE_j + \RR_{m+1+m_0;\,\delta}.
\end{equation}
Moreover by definition the $\EE^+_j$ do not depend on $\delta$ and the $\EE^-_j$ are profiles: Using \eqref{6E4b} we find that for any $j\in\N$
\begin{equation}
\label{6E16b}
   \delta^{-\frac12} \|\EE^-_j\|_{0,\Omega_-} + 
   \delta^{\frac12} \|\rot\EE^-_j\|_{0,\Omega_-} \leqslant C
   \|\WW_j\|_{\bH(\rot,\Sigma\times\R_+)} .
\end{equation}
Combining with \eqref{6E4c}, we obtain for any $j\in\N$
\begin{equation}
\label{6E17}
   \|\EE^+_j\|_{\bH(\rot,\Omega_+)} + 
   \delta^{-\frac12} \|\EE^-_j\|_{0,\Omega_-} + 
   \delta^{\frac12} \|\rot\EE^-_j\|_{0,\Omega_-} \leqslant C_j.
\end{equation}
We finally deduce the wanted estimate \eqref{6E8} from \eqref{6E15} to \eqref{6E17}.
\end{proof}

\begin{rem}
\label{6R1}
As a consequence of the works \cite{H-J-N07} and \cite{Peron09}, we find the existence of asymptotics of the form \eqref{6E4} when the interface $\Sigma$ is smooth, if the right hand side $\jj$ is smooth and has its support in the dielectric part $\Omega_+$. Moreover, estimate \eqref{6E7} is true for $m_{0}=1$, cf.\ \cite[Ch.\,7]{Peron09}.
\end{rem}

\begin{rem}
\label{6R2}
If the interface has conical points, or is polyhedral, many difficulties are encountered for an asymptotic analysis. We refer to \cite{MaghnoujiNicaise06} for an investigation of a scalar transmission problem with high contrast in polygonal domain.
\end{rem}

\end{document}